\documentclass[aos,seceqn,citesort,dvips]{arximspdf}

\doi{10.1214/09-AOS770}
\volume{38}
\issue{4}
\pubyear{2010}
\firstpage{2047}
\lastpage{2091}

\makeatletter

\newtheorem{fa}{Fact}
\newproclaim{Definition}{Definition}
\newproclaim{Remark}{Remark}
\newtheorem{lm}{Lemma}
\newtheorem{theorem}{Theorem}[section]
\newtheorem{cor}{Corollary}
\newtheorem{ex}{Example}

\makeatother

\begin{document}
\begin{frontmatter}

\title{Testing conditional independence using maximal nonlinear
conditional correlation\protect\thanksref{T1,T2,T3}}
\runtitle{Testing conditional independence}

\thankstext{T1}{A part of this work is done when
the author was visiting the Institute of Statistical Science at Academia
Sinica in Taiwan.}

\thankstext{T2}{Supported by the National Science Council of Taiwan
through Grants
NSC 95-2119-M-004-001- and NSC 97-2118-M-004-001-.}

\thankstext{T3}{Supported in part by Center for Service Innovation at
National Chengchi University.}

\begin{aug}
\author[A]{\fnms{Tzee-Ming} \snm{Huang}\corref{}\ead[label=e1]{tmhuang@nccu.edu.tw}}
\runauthor{T.-M. Huang}
\affiliation{National Chengchi University}
\address[A]{Department of Statistics\\
National Chengchi University\\
NO. 64, Sec. 2, ZhiNan Rd.\\
Taipei 11605, Taiwan\\
Republic of China\\
\printead{e1}} 
\end{aug}

\received{\smonth{9} \syear{2009}}
\revised{\smonth{11} \syear{2009}}

%
\begin{abstract}
In this paper, the maximal nonlinear conditional correlation of two
random vectors $X$ and $Y$ given another random vector $Z$, denoted by
$\rho_1(X,Y|Z)$, is defined as a measure of conditional association,
which satisfies certain desirable properties. When $Z$ is continuous, a
test for testing the conditional independence of $X$ and $Y$ given $Z$
is constructed based on the estimator of a weighted average of the form
$\sum_{k=1}^{n_Z} f_Z(z_k) \rho^2_1(X,Y|Z=z_k)$, where $f_Z$ is the
probability density function of $Z$ and the $z_k$'s are some points in
the range of $Z$. Under some conditions, it is shown that the test
statistic is asymptotically normal under conditional independence, and
the test is consistent.
\end{abstract}

%

\end{frontmatter}

\section{Introduction}\label{sec1}
In this paper, the problem of interest is testing the conditional
independence between two random vectors $X$ and $Y$
given a third random vector $Z$. The study of the problem of testing
conditional independence has a long history.
However, there are relatively few results on nonparametric tests when
the vectors $X$, $Y$ and $Z$ are continuous.
Some examples of such tests can be found in Su and White \cite
{MR2413488,MR2428851}, 
where they also proposed conditional independence
tests based on a weighted Hellinger distance between the conditional
densities or the difference between the
conditional characteristic functions.

As mentioned in Daudin \cite{MR601095}, 
$X$ and $Y$ are conditionally independent given $Z$ means that for
every $f(X,Z)$ and $g(Y,Z)$ such that
$Ef^2(X,Z)$ and $Eg^2(Y,Z)$ are finite
\[
E(f(X,Z)g(Y,Z)|Z) = E(f(X,Z)|Z) E(g(Y,Z)|Z).
\]
Thus, the problem of testing conditional independence,
as the problem of testing unconditional independence, is invariant when
one-to-one transforms are applied to the marginals $X$ and $Y$,
respectively. Various authors have taken this invariant
property into consideration when constructing conditional or
unconditional independence tests.
For example, Su and White \cite{MR2428851} 
used Hellinger distance in their test statistic for testing conditional
independence, so
that the test statistic is invariant. Dauxois and Nkiet \cite
{MR1647653} 
used measures of association to construct independence tests, and
the measures are invariant under the above transforms. In this paper,
to take invariance into account, the proposed test is
based
on the maximal nonlinear conditional correlation, which can be viewed
as a measure of conditional association and satisfies the
above invariance property.

To choose a reasonable measure of conditional association between $X$
and $Y$, the following properties are considered.
\begin{enumerate}[(P5)]
\item[(P1)] The measure can be defined for all types of random vectors,
including both discrete and continuous ones.
\item[(P2)] The measure is symmetric, that is, it remains the same when
$(X,Y)$ is replaced by $(Y,X)$.
\item[(P3)] The measure is invariant when one-to-one transforms are
applied to $X$ and~$Y$, respectively.
\item[(P4)] The measure is between 0 and 1.
\item[(P5)] The measure is 0 if and only if conditional independence holds.
\end{enumerate}
The above properties are adapted from some of the conditions for a good
measure of association proposed by
R\'{e}nyi \cite{MR0115203}. 
In \cite{MR0115203}, the conditional independence in (P5) is replaced
by the unconditional
independence.
Note
that the symmetric property (P2) is not always required. For instance,
Hsing et al. \cite{hsingcoefficient2005} 
proposed to use the coefficient of intrinsic dependence as a measure of
dependence, which does not satisfy (P2). Here, (P2) is considered.



Many measures of conditional association satisfying (P1)--(P5) can be
constructed. Dauxois and Nkiet \cite{MR1921388} 
showed
that a class of measures of association between two Hilbertian
subspaces can be obtained by properly combining the canonical
coefficients of the canonical analysis (CA) between the spaces. In
particular, take the two subspaces to be
$\tilde{H}_1 = \{ f(X,Z) - E(f(X,Z)|Z)\dvtx Ef^2(X,Z) < \infty\}$ and
$\tilde{H}_2 = \{ g(Y,Z) - E(g(Y,Z)|Z)\dvtx\break Eg^2(Y,Z) < \infty
\}$, 
then a class of measures of conditional association between $X$ and $Y$
given $Z$ satisfying properties (P1)--(P5) can
be obtained using the canonical coefficients. Denote the canonical
coefficients (arranged in descending order) by
$\tilde{\rho}_{i}(X,Y|Z)\dvtx i=1, 2, \ldots.$ When $X$ and $Y$ are not
functions of $Z$, the largest canonical
coefficient $\tilde{\rho}_{1}(X,Y|Z)$ is the maximal partial
correlation defined by Romanovi\v{c} \cite{MR0420757}, 
which is
\[
\sup_{f, g} \operatorname{corr}  \bigl( f(X,Z) - E(f(X,Z)|Z), g(Y,Z) -
E(g(Y,Z)|Z)  \bigr).
\]
%

Another approach to construct measures of conditional association is to
modify the CA between the spaces $H_1 = \{ f(X)
- Ef(X)\dvtx Ef^2(X) < \infty\}$ and $H_2 = \{ g(Y) - Eg(Y)\dvtx Eg^2(Y) <
\infty\}$ to obtain a conditional version of it.
That is, to find pairs of functions $(f_i, g_i)\dvtx i=0, 1, \ldots,$
such that for each $i$, $(f_i,g_i)$ maximizes
$E ( f(X,Z)g(Y,Z)|Z  )$ subject to
%
%
%
\begin{eqnarray}
\label{eq:n1}\quad
E(f^2(X,Z)|Z) I_{(0, \infty)} (E(f^2(X,Z)|Z))
&=& I_{(0, \infty)} (E(f^2(X,Z)|Z)),
\\
\label{eq:n2}
E(g^2(Y,Z)|Z) I_{(0, \infty)}
(E(g^2(Y,Z)|Z)) &=& I_{(0, \infty)} (E(g^2(Y,Z)|Z))
\end{eqnarray}
and
\[
E(f(X,Z)f_j(X,Z)|Z)
= 0 = E(g(Y,Z)g_j(Y,Z)|Z) \qquad\mbox{for } 0 \leq j < i.
\]
Here, $I_A$ denotes the indicator function on a set
$A$, that is, $I_A(x) = 1$ if $x \in A$ and $I_A(x) = 0$, otherwise. If
the above $(f_i, g_i)$'s exist, then one can define
$\rho_i(X,Y|Z) = E ( f_i(X,Z)g_i(Y,Z)|Z  )$ for each $i$ and
the $\rho_i(X,Y|Z)$'s can serve as a conditional
version of canonical coefficients. A measure of conditional association
satisfying (P1)--(P5)
can be obtained by taking a proper combination of the $\rho
_i(X,Y|Z)$'s, following the approach in \cite{MR1921388}. 
Examples of such combinations
include $\rho_1(X,Y|Z)$ and $1- \exp( - \sum_i \rho^2_i(X,Y|Z)
)$. The measure of conditional association used in this paper is $\rho
_1(X,Y|Z)$, which
will be called the maximal nonlinear conditional correlation of
two random vectors $X$ and $Y$ given $Z$ from now on.

In the above definition of $\rho_i(X,Y|Z)$'s, it is assumed that the
$(f_i, g_i)$'s exist. However, it is not clear what conditions can
guarantee the existence
of the $(f_i, g_i)$'s. To avoid the problem of finding
such conditions, a more general definition
for $\rho_1(X,Y|Z)$ is given in Section \ref{sec2}.
To construct a test based on $\rho_1(X,Y|Z)$, it is assumed that $Z$
has a Lebesgue probability density function $f_Z$.
An estimator of $\sum_{k} f_Z(z_k) \rho^2_1(X,Y|Z=z_k)$
is then used as the test statistic, where the $z_k$'s are some points
in the range of $Z$.
To study the asymptotic behavior of the test statistic under the
hypothesis that $X$ and $Y$ are conditionally independent given $Z$, we
follow the
approach in \cite{MR1647653} 
for finding the asymptotic distribution of a statistic for testing the
independence between $X$ and $Y$, which is
based on estimators of the canonical coefficients from the CA of $H_1$
and $H_2$. To make the
approach work for the conditional case,
some strong approximation results for kernel estimators of certain
conditional expectations are also established.

This paper is organized as follows. The new definition of $\rho
_1(X,Y|Z)$ is given in Section \ref{sec2}. Section \ref{sec3} deals with the estimation
of $\rho_1(X,Y|Z=z)$ and test construction. An example is in Section
\ref{sec:ex} and proofs are given in Section \ref{sec:proofs}.

\section{Maximal nonlinear conditional correlation}\label{sec2}
In this section, a more general definition of the maximal nonlinear
conditional correlation $\rho_1(X,Y|Z)$ will be given.
Note that in the definition of $\rho_i(X,Y|Z)$'s in Section \ref{sec1}, one can
take $f_0(X,Z)=1=g_0(Y,Z)$, which gives
that $\rho_0(X,Y|Z)=1$, and then
$\rho_1(X,Y|Z)$ can be defined as $E(f_1(X,Z)g_1(Y,Z)|Z)$ if there
exists $(f_1, g_1) \in S_0$ such that
\[
E(f(X,Z)g(Y,Z)|Z) \leq E(f_1(X,Z)g_1(Y,Z)|Z) \qquad\mbox{for every } (f,g)
\in S_0,
\]
where $S_0$ is the collection of pairs of functions $(f, g)$'s that
satisfy (\ref{eq:n1}), (\ref{eq:n2}) and
$E(f(X,Z)|Z) =0 = E(g(Y,Z)|Z)$.
Without assuming the existence of $(f_1, g_1)$, it is reasonable to
define $\rho_1(X, Y|Z)$ as
%
%
\begin{equation} \label{eq:sup1} \sup_{(f,g) \in S_0}
E(f(X,Z)g(Y,Z)|Z),
\end{equation}
if the supremum can be defined.

The above approach can be considered as a ``pointwise'' approach.
Indeed, when $Z$ takes values in a countable set $\mathcal{Z}$, for
each $z \in
\mathcal{Z}$, one may define \mbox{$\rho_1(X, Y|Z=z)$} as
%
%
\begin{equation} \label{eq:sup2}
\sup_{(f,g) \in S_0}E\bigl(f(X,z)g(Y,z)|Z=z\bigr),
\end{equation}
then the $\rho_1(X, Y|Z)$ defined using (\ref{eq:sup2}) is a measurable
function and can serve as the supremum in (\ref{eq:sup1}). However, if
$\mathcal{Z}$ is
uncountable, then it is not clear whether the $\rho_1(X, Y|Z)$ defined
using (\ref{eq:sup2}) is measurable. Therefore, we use the following
fact to define
the supremum in (\ref{eq:sup1}) so that it is well defined and is a
measurable function.
\begin{fa} \label{fa:1} There exists a sequence $\{ (\alpha_n, \beta_n)
\}$ in $S_0$ such that:
\begin{longlist}
\item The sequence $\{ E(\alpha_n(X,Z) \beta_n(Y,Z)|Z) \}$ is
nondecreasing, and
\item for every $(f,g) \in S_0$,
\[
E(f(X,Z)g(Y,Z)|Z) \leq\lim_{n\rightarrow\infty}
E(\alpha_n(X,Z) \beta_n(Y,Z)|Z).
\]
\end{longlist}
Furthermore, if \textup{(i)} and \textup{(ii)} hold for $\{ (\alpha_n , \beta_n) \} = \{
(\alpha_{n,1}, \beta_{n,1}) \}$
or $\{ (\alpha_{n,2}, \beta_{n,2}) \}$, where $\{ (\alpha_{n,1},
\beta
_{n,1}) \}$ and $\{ (\alpha_{n,2}, \beta_{n,2}) \}$
are sequences in $S_0$, then
%
%
\begin{equation} \label{eq:fact1-1}\qquad
\lim_{n \rightarrow\infty}
E (\alpha_{n,1} (X,Z)\beta_{n,1} (Y,Z)|Z) = \lim_{n \rightarrow
\infty}
E (\alpha_{n,2}(X,Z)\beta_{n,2}(Y,Z)|Z).
\end{equation}
\end{fa}

For the sake of brevity, from now on, some functions of $(X,Z)$ or
$(Y,Z)$ may
be expressed without the arguments $(X,Z)$ or $(Y,Z)$. For
distinguishing purpose, functions of $(X,Z)$ may
have names starting with only $\alpha$ or $f$, and functions of $(Y,Z)$
may have names starting with only
$\beta$ or $g$.
\begin{pf*}{Proof for Fact \ref{fa:1}}
We will first establish (\ref{eq:fact1-1}) if (i) and (ii) hold for $\{
(\alpha_n , \beta_n) \} = \{ (\alpha_{n,1},
\beta_{n,1}) \}$
or $\{ (\alpha_{n,2}, \beta_{n,2}) \}$. Note that for each $n$, from
(ii), we have that
\[
E (\alpha_{n,2}\beta_{n,2}|Z) \leq\lim_{n \rightarrow\infty} E
(\alpha_{n,1} \beta_{n,1}|Z)
\]
and
\[
E (\alpha_{n,1}\beta_{n,1} |Z) \leq\lim_{n \rightarrow\infty} E
(\alpha_{n,2}\beta_{n,2}|Z).
\]
Take the limits in these two inequalities as $n\rightarrow\infty$, and
we have (\ref{eq:fact1-1}).

It remains to find a sequence $\{ (\alpha_n, \beta_n) \}$ in $S_0$ that
satisfies (i) and (ii). Let
$\{ (\alpha_{n,0}, \beta_{n,0}) \}$ be a sequence in $S_0$ so that the
sequence $\{ E(\alpha_{n,0}\beta_{n,0}) \}$ is
nondecreasing and converges to $\sup_{(f,g) \in S_0} E (fg)$. We will
construct $\{ (\alpha_n, \beta_n) \}$ using
$\{ (\alpha_{n,0}, \beta_{n,0}) \}$ as follows. For $n=1$, define
$(\alpha_1,\beta_1) = (\alpha_{1,0}, \beta_{1,0})$.
For $n \geq2$, define
\begin{eqnarray*}
&&(\alpha_n(X,Z),\beta_n(Y,Z)) \\
&&\qquad =  \cases{(\alpha_{n,0}(X,Z), \beta_{n,0}(Y,Z)),
&\quad if $E (\alpha_{n,0}\beta_{n,0}|Z) > E(\alpha_{n-1}\beta
_{n-1}|Z)$; \cr
(\alpha_{n-1}(X,Z), \beta_{n-1}(Y,Z)), &\quad otherwise.}
\end{eqnarray*}
Then $\{ (\alpha_n,\beta_n) \}$ is a sequence in $S_0$ that satisfies
(i), and the sequence $\{ E\alpha_n\beta_n \}$
converges to $\sup_{(f,g) \in S_0} E (fg)$ since $E(\alpha_n\beta_n|Z)
\geq E(\alpha_{n,0}\beta_{n,0}|Z)$.
To see that $\{ (\alpha_n,\beta_n) \}$ also satisfies (ii), for
$(\alpha
, \beta)$ in $S_0$, define
\[
(\alpha_n^{*}, \beta_n^{*}) =  \cases{
(\alpha,\beta), &\quad if $E (\alpha\beta|Z) > \displaystyle\lim
_{n\rightarrow
\infty} E(\alpha_n \beta_n |Z)$; \cr
(\alpha_n, \beta_n), &\quad otherwise.}
\]
Then $\{ (\alpha_n^*,\beta_n^*) \}$ is a sequence in $S_0$ such that
%
%
\begin{equation} \label{eq:fa:1-1}
\lim_{n\rightarrow\infty} E(\alpha_n^{*}\beta_n^{*}|Z) = \max
\Bigl\{
E(\alpha\beta|Z), \lim_{n\rightarrow\infty}
E(\alpha_n \beta_n|Z)  \Bigr\}.
\end{equation}
From the monotone convergence theorem, we have
%
%
\begin{equation} \label{eq:fa:1-mct-1}
E\lim_{n\rightarrow\infty} E(\alpha_n^{*}\beta_n^{*}|Z) = \lim
_{n\rightarrow\infty} E(\alpha_n^{*}\beta_n^{*})
\end{equation}
and
%
%
\begin{equation} \label{eq:fa:1-mct-2} E\lim_{n\rightarrow\infty
}E(\alpha_n \beta_n|Z)
= \lim_{n\rightarrow\infty} E(\alpha_n \beta_n),
\end{equation}
so (\ref{eq:fa:1-1}) implies that
\[
\sup_{(f,g) \in S_0} E(fg) \geq\lim_{n\rightarrow\infty} E(\alpha
_n^{*}\beta_n^{*}) \geq
\lim_{n\rightarrow\infty} E(\alpha_n \beta_n) = \sup_{(f,g) \in S_0}
E(fg),
\]
which gives
%
%
\begin{equation}
\label{eq:fa:1-2} \lim_{n\rightarrow\infty} E(\alpha_n^{*}\beta
_n^{*}) = \lim_{n\rightarrow\infty}
E(\alpha_n \beta_n).
\end{equation}
If $E(\alpha\beta|Z) > \lim_{n\rightarrow\infty} E(\alpha_n \beta
_n |Z)$
with positive probability, then (\ref{eq:fa:1-1}),
(\ref{eq:fa:1-mct-1}) and (\ref{eq:fa:1-mct-2}) together
implies that $\lim_{n\rightarrow\infty} E(\alpha_n^{*}\beta_n^{*}) >
\lim_{n\rightarrow\infty}E(\alpha_n \beta_n)$,
which contradicts (\ref{eq:fa:1-2}). Thus, (ii) holds. The proof of
Fact \ref{fa:1} is complete.
\end{pf*}

With Fact \ref{fa:1}, the maximal nonlinear conditional correlation
$\rho_1(X,Y|Z)$ can be redefined as follows.
\begin{Definition}\label{defi1}
$\rho_1(X,Y|Z) = \sup_{(f,g) \in S_0}
E(f(X,Z)g(Y,Z)|Z)$, which is defined as
$\lim_{n\rightarrow\infty} E(\alpha_n(X,Z) \beta_n(Y,Z)|Z)$, where
$\{
(\alpha_n, \beta_n) \}$ is a sequence in
$S_0$ that satisfies (i) and (ii) in Fact \ref{fa:1}.
\end{Definition}

Below are some remarks for the $\rho_1(X,Y|Z)$.
\begin{enumerate}
\item If there exists $(f_1, g_1)$ in $S_0$ such that $E(f_1g_1|Z) \geq
E(fg|Z)$ for all $(f,g) \in S_0$, then
$\rho_1(X,Y|Z) = E(f_1g_1|Z)$ using Definition \ref{defi1}. To see this, let $\{
(\alpha_n, \beta_n) \}$ be a sequence in
$S_0$ that
satisfies (i) and (ii) in Fact \ref{fa:1}. Then $\rho_1(X,Y|Z) = \lim
_{n\rightarrow\infty} E(\alpha_n\beta_n|Z)$, so $E(f_1g_1|Z) \leq
\rho
_1(X,Y|Z)$ by (ii). Also,
$E(f_1g_1|Z) \geq E(\alpha_n\beta_n|Z)$ for every $n$, so $E(f_1g_1|Z)
\geq\rho_1(X,Y|Z)$. Therefore, $\rho_1(X,Y|Z) = E(f_1g_1|Z)$ and
Definition \ref{defi1} can be viewed as a generalized version
of the definition of $\rho_1(X,Y|Z)$ given in Section \ref{sec1}.

\item$\rho_1(X,Y|Z)$ satisfies properties (P1)--(P5).

\item When $X$ is a function of $Y$ and $Z$ or $Y$ is a function of $X$
and $Z$, it is not necessary that
$\rho_1(X,Y|Z) = 1$. For instance, suppose that $X$ and $Z$ are
independent standard normal random variables
and $Y=XI_{(0,\infty)}(Z)$, then $\rho_1(X,Y|Z) = I_{(0,\infty)}(Z)$.

\item Let $\rho_1(X,Y)$ be the largest canonical
coefficient from the CA between $H_1 = \{ f(X)
- Ef(X)\dvtx Ef^2(X) < \infty\}$ and $H_2 = \{ g(Y) - Eg(Y)\dvtx Eg^2(Y) <
\infty\}$. Then $\rho_1(X,Y|Z) = \rho_1(X,Y)$
if $(X, Y)$ and $Z$ are independent.

\item Let $\rho_1(X,Y)$ be as defined in item 4. It is stated
in \cite{MR1647653} 
that when the joint distribution of $X$ and $Y$ is bivariate normal
\[
N \left(
\pmatrix{ 0 \cr 0 },
\pmatrix{ 1 & \rho\cr \rho& 1}
\right),
\]
$\rho_1(X,Y) = |\rho|$. This result implies that, when the joint
distribution for $X$, $Y$ and $Z$ is multivariate normal and $X$ and
$Y$ are both univariate,
\begin{eqnarray*}
\rho_1(X,Y|Z)
& = &  \biggl| \frac{ E( (X- E(X|Z))(Y-E(Y|Z))|Z ) }{  (
E(X-E(X|Z))^2|Z)^{1/2}
 ( E(Y-E(Y|Z))^2|Z)^{1/2}}  \biggr| \\
& = &  \biggl| \frac{ E(X-E(X|Z))(Y-E(Y|Z)) }{  (
E(X-E(X|Z))^2  )^{1/2}  ( E(Y-E(Y|Z))^2  )^{1/2}}
\biggr|,
\end{eqnarray*}
which also equals the absolute value of the usual partial correlation
coefficient.
\end{enumerate}

\section{A test of conditional independence}\label{sec3}
Testing conditional independence is equivalent to testing $H_0\dvtx \rho
_1(X,Y|Z)=0$, which involves testing $H_{0,z}\dvtx\break \rho_1(X,Y| Z=z)=0$ for
different $z$'s in the range of $Z$. Let $\mathcal{Z}$ be the range of
$Z$. In
this section, an estimator $\hat{\rho}(z)$ is proposed for estimating
$\rho_1(X,Y|Z=z)$ for each $z \in\mathcal{Z}$, and for distinct
points $z_1,
\ldots, z_{n_Z}$ in $\mathcal{Z}$, the asymptotic joint
distribution of $\hat
{\rho}(z_1), \ldots, \hat{\rho}(z_{n_Z})$ under $H_0$ is
derived to
construct a test for testing $H_0$.

\subsection{Estimation of $\rho_1(X,Y|Z=z)$} \label{section:estimator}
To estimate
\[
\rho_1(X,Y|Z) = \sup_{ (f, g) \in S_0} E(fg|Z)
\]
for $(f,g) \in S_0$, $f$ and $g$ are approximated using basis functions.
Suppose that there exist $\Lambda_1$, $\Lambda_2$ and $\Lambda_3$:
subsets of the set of all positive integers and three sets
of functions $\{ \phi_{p, i}\dvtx 1 \leq i \leq p, p \in\Lambda_1 \}$,
$\{
\psi_{q, j}\dvtx
1 \leq j \leq q, q \in\Lambda_2 \}$ and $\{ \theta_{r, k} \dvtx
1 \leq k \leq r, k \in\Lambda_3 \}$ such that for $\alpha(X,Z)$ and
$\beta(Y,Z)$ with finite second moments,
%
%
\begin{equation}\quad
\label{eq:app1} \lim_{p,r\rightarrow\infty} \inf_{ a(i,k) }
E \biggl( \alpha(X,Z) - \sum_{ 1\leq i \leq p, 1 \leq k \leq r }
a(i,k) \phi_{p,i}(X) \theta_{r,k}(Z)  \biggr)^2 = 0
\end{equation}
and
%
%
\begin{equation} \label{eq:app2}\quad
\lim_{q,r\rightarrow\infty} \inf_{ b(j,k) } E \biggl( \beta(Y,Z) -
\sum
_{1\leq j \leq q, 1 \leq k \leq r}
b(j,k) \psi_{q,j}(Y) \theta_{r,k}(Z)  \biggr)^2 = 0.
\end{equation}
Also, suppose that for each $(p,q)$, there
exist coefficients $a_{p,0,i}$'s and $b_{q,0,j}$'s such that
%
%
\begin{equation} \label{eq:basis}
\sum_{1\leq i \leq p}
a_{p,0, i} \phi_{p,i}(x) =1 = \sum_{1\leq j \leq q} b_{q, 0,j} \psi
_{q,j}(y)
\end{equation}
for every $x$ in the range of $X$ and every $y$ in the range of $Y$.

Let $S_1$ be the collection of all $(f,g)$'s with finite second moments
and let $S_{1,p,q}$ be the collection of all $(f,g)$'s in $S_1$
such that $f(X,Z) =\break \sum_{i=1}^{p} a_{p,i} (Z)
\phi_{p, i} (X)$ for some $a_{p,i} (Z)$'s, and $g(Y,Z) = \sum_{j=1}^{q}
b_{q,j}(Z) \psi_{q, j}(Y)$ for
some
$b_{q,j}(Z)$'s. Then (\ref{eq:app1}) and (\ref{eq:app2}) together
imply that $S_1$ can be approximated by $S_{1,p,q}$ for large $p$ and~$q$.
Since $S_0 \subset S_1$, $S_0$ can be approximated by $S_{1,p,q}$ as
well. With the additional condition (\ref{eq:basis}),
$S_0$ can be easily approximated using the subspace $S_{0,p,q} = S_0
\cap S_{1, p,q}$. Note that (\ref{eq:app1}), (\ref{eq:app2}) and
(\ref{eq:basis}) hold for certain basis functions, for example, the
tensor product splines in \cite{schumakerspline1981}. 

Assuming (\ref{eq:app1}), (\ref{eq:app2}) and (\ref{eq:basis}), it is
reasonable to define
\[
\sup_{ (f, g) \in S_{0,p,q}} E(fg|Z)
\]
and
use it to approximate $\rho_1(X,Y|Z)$. To define $\sup_{ (f, g) \in
S_{0,p,q}} E(fg|Z)$, one may follow the same approach for defining
$\sup_{ (f, g) \in S_0} E(fg|Z)$, or simply note that there exists
$(f_1, g_1) \in S_{0, p,q}$ such that
%
%
\begin{equation} \label{eq:max} E(f_1g_1|Z) \geq E(fg|Z) \qquad\mbox{for all
$(f, g) \in S_{0, p,q}$}
\end{equation}
and define $\sup_{ (f, g) \in S_{0,p,q}} E(fg|Z) = E(f_1g_1|Z)$. The
pair $(f_1, g_1)$ can be obtained as follows. Let
\begin{eqnarray*}
\Sigma_{\phi, p}(Z) &=&  \bigl( E(\phi_{p, i}(X) \phi_{p, j}(X)|Z)-
E(\phi
_{p, i}(X)|Z) E(\phi_{p, j}(X)|Z)  \bigr)_{p \times p},
\\
\Sigma_{\psi, q} (Z) &=&  \bigl( E(\psi_{q, i}(Y) \psi_{q, j}(Y)|Z) -
E(\psi
_{q, i}(Y)|Z) E(\psi_{q, j}(Y)|Z)  \bigr)_{q \times q}
\end{eqnarray*}
and
\[
\Sigma_{\phi, \psi, p,q}(Z) =  \bigl( E(\phi_{p, i}(X) \psi_{q,
j}(Y)|Z)- E(\phi_{p, i}(X)|Z) E(\psi_{q, j}(Y)|Z)  \bigr)_{p \times q}.
\]
Consider the following two cases:
\begin{longlist}
\item $\Sigma_{\phi, p}(Z)$ and $\Sigma_{\psi, q} (Z)$ are not
zero matrices, and
\item at least one of $\Sigma_{\phi, p}(Z)$ and $\Sigma_{\psi, q}
(Z)$ is a zero matrix.
\end{longlist}
In case (i), let $a_1 =(a_{1,1}(Z), \ldots, a_{1,p}(Z))^T$ and $b_1=
(b_{1,1}(Z), \ldots, b_{1,q}(Z))^T$ be such that $(a_1, b_1)$ is the
pair of $(a,b)$ that
maximizes
\[
a^T \Sigma_{\phi, \psi, p,q}(Z) b
\]
subject to
\[
a^T\Sigma_{\phi, p}(Z) a =1= b^T \Sigma_{\psi, q} (Z) b,
\]
and then take
\[
f_1(X,Z) = \sum_{i=1}^p a_{1,i}(Z) \bigl( \phi_{p,i}(X) - E(\phi_{p,i}(X)|Z)
\bigr)
\]
and
\[
g_1(Y,Z) = \sum_{j=1}^q b_{1,j}(Z) \bigl( \psi_{q,j}(Y) - E(\psi_{q,j}(Y)|Z)
\bigr).
\]
In case (ii), take $f_1(X,Z) = 0 = g_1(Y,Z)$. Then $(f_1, g_1) \in
S_{0, p, q}$ and (\ref{eq:max}) holds. Denote $\sup_{ (f, g) \in
S_{0,p,q}} E(fg|Z)$ by
$\rho_{p, q}(Z)$.


The following fact states that $\rho_1(X,Y|Z)$
can be reasonably approximated by $\rho_{p, q}(Z)$ if $p$ and $q$ are large.
\begin{fa} \label{fact:app}
Suppose that (\ref{eq:app1}), (\ref{eq:app2}) and (\ref{eq:basis}) hold
and $\{ p_n \}$ and $\{ q_n \}$ are
sequences of positive integers that tend to $\infty$ as $n \rightarrow
\infty
$. Then
\[
\lim_{n\rightarrow\infty} E\bigl(|\rho_1(X,Y|Z) - \rho_{p_n, q_n}(Z)|\bigr)
= 0.
\]
\end{fa}
\begin{pf}
Since $\rho_1(X,Y|Z) \geq\rho_{p_n,
q_n}(Z)$ for every $n$, Fact \ref{fact:app} holds if for every
$\varepsilon> 0$, there exists
$N_0$ such that for $n \geq N_0$,
%
%
\begin{equation} \label{eq:rho-d1} \rho_1(X,Y|Z) \leq\rho_{p_n,
q_n}(Z) +\Delta_1
\end{equation}
for some $\Delta_1$ such that $E|\Delta_1| < \varepsilon$. To find such
a $\Delta_1$, we will first look for a pair $(f_m, g_m) \in S_0$
such that $E(f_mg_m|Z) \approx\rho_1(X,Y|Z)$, and then find $(f^*_n,
g^*_n) \in S_{0,p_n,q_n}$ such that $(f^*_n, g^*_n) \approx(f_m, g_m)$.
Take
%
%
\begin{equation} \label{eq:d1} \Delta_1 = E(f_mg_m|Z) -E(f^*_ng^*_n|Z)
+ \rho_1(X,Y|Z) - E(f_mg_m|Z),
\end{equation}
then (\ref{eq:rho-d1}) holds
and $E|\Delta_1|$ can be made small if $m$ and $n$ are large enough.

To find $(f_m, g_m) \in S_0$ such that $E(f_mg_m|Z) \approx\rho
_1(X,Y|Z)$, let $\{ (f_n, g_n) \}_{n=1}^{\infty}$ be a sequence in
$S_0$ such that
$\{ E(f_ng_n|Z) \}$ is an increasing sequence and $\lim_{n\rightarrow
\infty}
E(f_ng_n|Z) = \rho_1(X,Y|Z)$.
Let $\Delta_{2,n}
= \rho_1(X,Y|Z) - E(f_ng_n|Z)$, then $\lim_{n\rightarrow\infty}
E|\Delta
_{2,n}|=0$, which
implies that for every $\delta>0$, there exists $m$ such that
%
%
\begin{equation} \label{eq:d2}
E|\Delta_{2,m}| < \delta.
\end{equation}

To find $(f^*_n, g^*_n) \in S_{0,p_n,q_n}$ such that $(f^*_n, g^*_n)
\approx(f_m, g_m)$, note that it follows from (\ref{eq:app1})
and (\ref{eq:app2}) that for $n \geq N_0$, there exists some
$(f_{n,1}, g_{n,1}) \in S_{1,p_n,q_n}$ such that
%
%
\begin{equation} \label{eq:app1-1} \sqrt{E(f_m - f_{n,1})^2} < \delta
\quad\mbox{and}\quad \sqrt{E(g_m - g_{n,1})^2} < \delta.
\end{equation}
Let $f_{n,2}(X,Z) = f_{n,1}(X,Z) - E (f_{n,1}|Z)$, $g_{n,2}(Y,Z) =
g_{n,1}(Y,Z) - E (g_{n,1}|Z)$,
\[
f_n^*(X,Z)
= \frac{f_{n,2}(X,Z)}{\sqrt{ E(f_{n,2}^2|Z) }} I_{(0, \infty)} (
E(f_{n,2}^2|Z))
\]
and
\[
g_n^*(Y,Z) = \frac{g_{n,2}(Y,Z) }{ \sqrt{E(g_{n,2}^2|Z) } } I_{(0,
\infty)} (E(g_{n,2}^2|Z)),
\]
then it follows from (\ref{eq:basis}) that
$(f_n^*, g_n^*) \in S_{0,p_n,q_n}$. To see that $(f^*_n, g^*_n) \approx
(f_m, g_m)$,
let $\Delta_3 = f_m -f_n^*$ and $\Delta_4 = g_m- g_n^*$, then it can be
shown that
%
%
\begin{equation} \label{eq:delta3}
E\Delta^2_3 \leq16 \delta^2 + 8 \delta
\end{equation}
and
%
%
\begin{equation} \label{eq:delta4}
E\Delta^2_4 \leq16 \delta^2 + 8 \delta.
\end{equation}
Below we will verify (\ref{eq:delta3}) only since the verification for
(\ref{eq:delta4}) is similar.
Write
$\Delta_3= f_m -f_{n,2} + f_{n,2} - f_n^*$, then by (\ref{eq:app1-1}),
%
%
\begin{equation} \label{eq:d31} E(f_m -f_{n,2} )^2 \leq4 \delta^2
\end{equation}
since
$E(f_m -f_{n,2} )^2 \leq2 (E(f_m - f_{n,1})^2 + E(f_{n,1} -
f_{n,2})^2 )$ and $(f_{n,1} - f_{n,2})^2 = ( E ((f_m - f_{n,1})|Z) )^2
\leq E( (f_m - f_{n,1})^2|Z)$.
Also,
\begin{eqnarray*}
E \bigl( (f_n^* - f_{n,2})^2|Z  \bigr) & = &  \bigl( 1-\sqrt{E
(f_{n,2}^2|Z)}  \bigr)^2 I_{(0, \infty)} ( E (f_{n,2}^2|Z) )\\
& \leq&  | 1 - E (f_{n,2}^2|Z)  | \\
& = &  \bigl| E \bigl( (f_m - f_{n,2})^2|Z  \bigr) - 2 E \bigl(
f_m(f_m -
f_{n,2}) |Z  \bigr)  \bigr| \\
& \leq& E \bigl( (f_m - f_{n,2})^2|Z  \bigr) + 2 \sqrt{E \bigl( (f_m -
f_{n,2})^2|Z  \bigr)},
\end{eqnarray*}
so
%
%
\begin{equation} \label{eq:d32}\qquad
E(f_{n,2} - f_n^*)^2 \leq E(f_m- f_{n,2})^2 + 2 \sqrt{ E(f_m -
f_{n,2})^2} \stackrel{\fontsize{8.36}{10.36}\selectfont{\mbox{(\ref{eq:d31})}}}{\leq} 4\delta^2 + 4\delta.
\end{equation}
Therefore, (\ref{eq:delta3}) follows from (\ref{eq:d31}), (\ref
{eq:d32}) and the inequality $E\Delta^2_3 \leq2 ( E(f_m -f_{n,2} )^2
+ E(f_{n,2} - f_n^*)^2 )$.

Finally, the $\Delta_1$ in (\ref{eq:d1}) is $E(f_n^*\Delta_4|Z) +
E(g_n^*\Delta_3|Z)
+ E(\Delta_3\Delta_4|Z) + \Delta_{2,m}$, so it follows from (\ref
{eq:delta3}),
(\ref{eq:delta4}), (\ref{eq:d2}) and the Cauchy inequality that
\[
E|\Delta_1| \leq3 \sqrt{16 \delta^2 + 8 \delta} + \delta.
\]
For $\varepsilon>0$, one can choose $\delta$ so that $3 \sqrt{16
\delta
^2 + 8 \delta} + \delta< \varepsilon$, then $E|\Delta_1| <
\varepsilon
$ as required.
The proof of Fact \ref{fact:app} is complete.
\end{pf}

Based on Fact \ref{fact:app}, it is reasonable to estimate $\rho
_1(X,Y|Z)$ using an estimator for
$\rho_{p, q}(Z)$, where $p$ and $q$ are large. To estimate $\rho_{p, q}(Z)$,
the following assumption is made:
\begin{enumerate}[(A1)]
\item[(A1)] There exists a version of the conditional distribution of
$(X,Y)$ given $Z$ such that
for every bounded function $g(X,Y)$, $E(g(X,Y)|Z)$ calculated using
that version is a continuous function of $Z$.
%
\end{enumerate}
From now on, we will use the version of conditional distribution in
(A1) to obtain $E(g(X,Y)|Z=z)$ for every bounded $g$
and every $z$ in the range of $Z$. It for each $(p, q)$, $1 \leq i \leq p$, $1 \leq j \leq q$, $| \phi
_{p, i} | \leq1$ and $|\psi_{q, j} | \leq1$, then each element in $\Sigma
_{\phi, p}(z)$, $\Sigma_{\psi, q}(z)$ and
$\Sigma_{\phi, \psi, p, q}(z)$ is a continuous function of~$z$, and
$\rho_{p, q}(z)$ is $\max_{a,b} a^T \Sigma_{\phi, \psi, p, q}(z) b$,
where the maximum is taken over all vectors $a$ and $b$ such that
\[
a^T \Sigma_{\phi, p}(z) a = 1 = b^T \Sigma_{\psi, q}(z) b.
\]

To estimate $\rho_{p, q}(z)$, we consider the estimator
\[
\hat{\rho}_{p,q}(z) = \max_{a,b} a^T \hat{\Sigma}_{\phi, \psi, p, q}
(z)b,
\]
where the maximum is taken over all vectors $a$ and $b$ such that
\[
a^T \hat{\Sigma}_{\phi, p}(z) a = 1 = b^T
\hat{\Sigma}_{\psi, q}(z) b,
\]
and $\hat{\Sigma}_{\phi, p}(z)$, $\hat{\Sigma}_{\phi, \psi, p, q}(z)$
and $\hat{\Sigma}_{\psi, q}(z)$ are obtained by replacing the
conditional expectations in $\Sigma_{\phi, p}(z)$,
$\Sigma_{\phi, \psi, p, q}(z)$ and $\Sigma_{\psi, q}(z)$ by their
kernel estimators. Specifically, each element in $\Sigma_{\phi, p}(z)$,
$\Sigma_{\phi, \psi, p, q}(z)$ and $\Sigma_{\psi, q}(z)$ is of the form
$E(UV|Z=z) - (E(U|Z=z))(E(V|Z=z))$, where $U$ and $V$ are
functions of $X$ or $Y$, so each of $E(UV|Z=z)$, $E(U|Z=z)$ and
$E(V|Z=z)$ is of the form $E(g(X,Y)|Z=z)$, which is estimated by
%
%
\begin{equation} \label{eq:kernel}
\hat{E} \bigl(g(X,Y)|Z=z\bigr) \stackrel{\mathrm{def}}{=} \frac{ \sum_{i=1}^{n}
g(X_i, Y_i) k_h(z-Z_i)}{\sum_{i=1}^{n} k_h(z-Z_i)},
\end{equation}
where $k_h(z) = h^{-d} k_0(z/h)$ and $k_0$ is a kernel function on
$R^d$ satisfying certain conditions which will be
specified later.
For each $z\in\mathcal{Z}$, to make $\hat{\rho}_{p, q}(z)$ a
reasonable estimator
for $\rho_1(X,Y|Z=z)$, we will take $p=p_n$, $q=q_n$ and $h=h_n$, where
$p_n \rightarrow\infty$, $q_n \rightarrow\infty$ and
$h_n \rightarrow0$ as $n \rightarrow\infty$. The estimator $\hat
{\rho}_{p_n,
q_n}(z)$ will be abbreviated as $\hat{\rho}(z)$ for each $z \in
\mathcal{Z}$.

The estimator $\hat{\rho}(z)$ can be expressed in a different form that
is easier to analyze. Let $X_*$ and $Y_*$ be
random vectors of length $p_n$ and $q_n$, respectively, such that given
the data $(X_1, Y_1, Z_1), \ldots, (X_n, Y_n, Z_n)$,
\[
(X_*^T,Y_*^T)= ( \phi_{p_n, 1}(X_\ell), \ldots, \phi_{p_n,
p_n}(X_{\ell
}), \psi_{q_n, 1}(Y_\ell), \ldots, \psi_{q_n, q_n}(Y_{\ell}) )
\]
with probability $k_h(z-Z_{\ell})/\sum_{i=1}^{n} k_h(z-Z_i)$ for $1
\leq\ell\leq n$. Then
$\hat{ \Sigma}_{\phi, \psi, p,q}(z) = EX_*Y_*^T - EX_* EY_*^T$,
$\hat{
\Sigma}_{\phi, p}(z) = EX_*X_*^T - EX_* EX_*^T$
and $\hat{ \Sigma}_{\psi, q}(z) = EY_*Y_*^T - EY_* EY_*^T$, where the
expectations are conditional expectations given the
data. Therefore, the estimator $\hat{\rho}(z)$ is the largest canonical
coefficient from the centered canonical analysis between
$X_*$ and $Y_*$.
Note that it follows from (\ref{eq:basis}) that
%
%
\begin{equation} \label{eq:basis-1}
a_{n,*}^TX_* =1 = b_{n,*}^T Y_*,
\end{equation}
where
\[
a_{n,*} = (a_{p_n,0, 1}, \ldots, a_{p_n,0, p_n})^T \quad\mbox{and}\quad b_{n,*}
= (b_{q_n,0, 1}, \ldots, b_{q_n,0, q_n})^T,
\]
so $\hat{\rho}(z)$ can also be obtained from the noncentered canonical
analysis between $X_*$ and $Y_*$.
Let
\begin{eqnarray*}
V_{1,1}(z) &=&  \bigl( E\bigl(\phi_{p_n, i}(X)\phi_{p_n, j} (X)|Z=z\bigr)  \bigr)_{p
_n \times p_n},
\\
V_{1,2}(z) &=&
\bigl( E\bigl(\phi_{p_n, i} (X)\psi_{q_n, j} (Y)|Z=z\bigr)  \bigr)_{p_n
\times
q_n},
\\
V_{2,2}(z)
&=&  \bigl( E\bigl(\psi_{q_n, i} (Y)\psi_{q_n, j}(Y)|Z=z\bigr)  \bigr)_{q_n
\times
q_n} \quad\mbox{and}\quad V_{2,1} (z) = V_{1,2}(z)^T
\end{eqnarray*}
for $1 \leq i, j \leq2$, let $\hat{V}_{i,j}(z)$ be the estimator of
$V_{i,j}(z)$ obtained by replacing the conditional expectations in
$V_{i,j}(z)$ by their kernel estimators as in (\ref{eq:kernel}). Then
$\hat{V}_{1,1}(z) = EX_*X_*^T$, $\hat{V}_{1,2}(z) = EX_*Y_*^T$,
$\hat{V}_{2,2}(z) = EY_*Y_*^T$, so $\hat{\rho}(z)$ is the square root
of the largest eigenvalue of the matrix
\[
\hat{V}_{1,2}(z) \hat{V}_{2,2}^{-1}(z) \hat{V}_{2,1}(z) \hat
{V}_{1,1}(z)^{-1} - \hat{V}_{1,1}(z) a_{n,*}a_{n,*}^T.
\]
Also,
$\rho_{p_n,q_n}(z)$ is the square root of the largest eigenvalue of
the matrix
\[
V_{1,2}(z) V_{2,2}^{-1}(z) V_{2,1}(z) V_{1,1}(z)^{-1} - V_{1,1}(z)
a_{n,*}a_{n,*}^T.
\]

To simplify the above matrix expressions, some notation is introduced
as follows.
For a $(p_n+q_n) \times(p_n + q_n)$ matrix $U$, express $U$ as
\[
\pmatrix{ U_{1,1} & U_{1,2} \cr U_{2,1} & U_{2,2}},
\]
where the dimension of $U_{1,1}$ is $p_n\times p_n$.
For $1 \leq i, j \leq2$, let $g_{i,j}$ be the mapping that maps $U$ to
$U_{i,j}$.
For a $p_n \times1$ vector $a$ and a $(p_n+q_n) \times(p_n + q_n)$
matrix $U$, define
\[
g(U, a) = g_{1,2}(U) g_{2,2}(U)^{-1} g_{2,1}(U) g_{1,1}(U)^{-1} -
g_{1,1}(U) a a^T,
\]
if $g_{2,2}(U)$ and $g_{1,1}(U)$ are invertible. Let
\[
V(z) =  \pmatrix{V_{1,1}(z) & V_{1,2}(z) \cr V_{2,1}(z) & V_{2,2}(z)}
\]
and
\[
\hat{V}(z) =
\pmatrix{\hat{V}_{1,1}(z) & \hat{V}_{1,2}(z) \cr \hat
{V}_{2,1}(z) & \hat{V}_{2,2}(z)},
\]
then $\hat{\rho}(z)$ is the square root of the largest eigenvalue of
$g(\hat{V}(z), a_{n,*})$
and $\rho_{p_n,q_n}(z)$ is the square root of the largest eigenvalue of
$g(V(z), a_{n,*})$.

The matrix $g(\hat{V}(z), a_{n,*})$ can be replaced by a different
matrix if basis change is performed. That is, suppose that
\[
\phi= (\phi_{p_n, 1}, \ldots, \phi_{p_n, p_n})^T \quad\mbox{and}\quad \psi=
(\psi_{q_n, 1}, \ldots, \psi_{q_n, q_n})^T
\]
are replaced by
$\phi^* = P_1 \phi$ and $\psi^* = Q_1 \psi$, respectively, and $\hat
{V}(z)$ becomes $\hat{V}^*(z)$ after such a change is made. Then
$\hat{\rho}(z)$ is also the square root of the largest eigenvalue of
the matrix $g(\hat{V}^*(z), \alpha^*)$, where
$\alpha^* =(P_1^{-1})^T a_{n,*}$ is a vector such that $(\alpha^*)^T
\phi^* = 1$. To make the expression for
$g(V^*(z), \alpha^*)$ simple, the matrices $P_1$ and $Q_1$ are chosen
so that
%
%
\begin{equation} \label{eq:new_basis_1}
\phi^*_1=1=\psi^*_1,
\end{equation}
$g_{1,1}(V^*(z))$ and $g_{2,2}(V^*(z))$ are identity matrices, and for
$1 \leq i \leq p_n$ and $1 \leq j \leq q_n$,
%
%
\begin{equation} \label{eq:new_basis_3}
E \bigl( \phi^*_i (X) \psi^*_j
(Y)|Z=z\bigr) = \delta_{i,j} \sqrt{\lambda_i},
\end{equation}
where $\phi^*_i$ and $\psi_j^*$ denote the $i$th element in $\phi^*$
and the $j$th element in $\psi^*$,
respectively, $\delta_{i,j}$ denotes the Kronecker symbol and the
$\lambda_i$'s are the eigenvalues of
$g(V^*(z), \alpha^*)$. Note that $(\alpha^*)^T = (1, 0, \ldots, 0)$
with the above choice of $P_1$ and $Q_1$.

%

\subsection{Asymptotic properties and a test of conditional
independence} \label{sec:test}
In this section, we will give asymptotic properties of the estimators
$\hat{\rho} (z_k)\dvtx 1 \leq k \leq n_Z$, where the $z_k$'s are distinct
points in $\mathcal{Z}$. First, we will establish the consistency of the
estimators, which relies on the fact that for each $k$, the two
matrices $g(\hat{V}^*(z_k), \alpha^*)$ and $g(V^*(z_k), \alpha^*)$ are
close, and their largest eigenvalues are $\hat{\rho}^2 (z_k)$ and
$\rho
^2_{p_n, q_n}(z_k)$. The difference between $g(\hat{V}^*(z_k), \alpha
^*)$ and $g(V^*(z_k), \alpha^*)$ depends on the difference of $\hat
{V}^*(z_k)$ and $V^*(z_k)$, and the difference between some conditional
expectation $E(g(X, Y, Z)|Z=z)$ and its kernel estimator
$\hat{E}(g(X, Y, Z)|Z=z) = \sum_{i=1}^n w_{0,i}(z) g(X_i,Y_i,z)/ \sum
_{i=1}^n w_{0,i}(z)$, where $w_{0,i}(z) = k_0(h_n^{-1}(z-Z_i))$. To
make it easier to derive the asymptotic properties of $\hat{E}(g(X, Y,
Z)|Z=z)$, some regularity conditions on the distribution of $(X,Y,Z)$
are imposed as follows.
\begin{enumerate}[(R3)]
\item[(R1)] There exists a $\sigma$-finite measure $\mu$ such that for
every $z \in\mathcal{Z}$, the conditional distribution of $(X,Y)$
given $Z=z$
has a p.d.f. $f(\cdot|z)$ with respect to $\mu$. Also, $Z$ has a Lebesgue
p.d.f. $f_Z$, and $f(x,y|z)$ and $f_Z(z)$ are twice differentiable with
respect to $z$.

\item[(R2)] There exists a function $h$ on $\mathcal{X}\times
\mathcal{Y}$ such that
\begin{eqnarray*}
&& \sup_{z\in\mathcal{Z}} \max \biggl(  | f(x,y|z)  |,
\max_{1 \leq
i \leq d} \biggl| \frac{\partial}{\partial z_i} f(x,y|z)  \biggr|,
\max_{1
\leq i,j \leq d} \biggl| \frac{\partial^2}{\partial z_i\, \partial z_j}
f(x,y|z)  \biggr|  \biggr) \\
&&\qquad \leq h(x,y)
\end{eqnarray*}
and $\int h(x,y) \,d\mu(x,y) < \infty$.

\item[(R3)] There exist constants $c_0$ and $c_1$ such that
\[
\sup_{z \in\mathcal{Z}} \max \biggl( | f_Z(z) |, \max_{1 \leq i
\leq d}  \biggl|
\frac{\partial}{\partial z_i} f_Z(z)  \biggr|, \max_{1 \leq i,j \leq
d} \biggl| \frac{\partial^2}{\partial z_i\, \partial z_j} f_Z(z)  \biggr|
\biggr) \leq c_0
\]
and $1/f_Z(z) \leq c_1$ for $z \in\mathcal{Z}$.
\end{enumerate}
Note that (R2) implies condition (A1) in Section
\ref{section:estimator}. For the kernel function $k_0$, conditions (K1)
and (K2) are assumed. The notation \mbox{$\| \cdot\|$} denotes the Euclidean
norm for a vector or the Frobenius norm for a matrix.
\begin{enumerate}[(K2)]
\item[(K1)] $k_0 \geq0$, $\sup_u k_0(u) < \infty$, $\int k_0(u) \,du
=1$, $\int u k_0(u) \,du =0$ and $\sigma_0^2 =\break \int\| u\|^2 k_0(u) \,du <
\infty$.

\item[(K2)] There exists positive constants $\gamma_2$ and $\gamma_3$
that does not depend on $d$ such that
\[
k_0(a) \leq(\gamma_2)^d e^{-\gamma_3 \| a \|^2} \qquad\mbox{for every $a
\in R^d$}.
\]
\end{enumerate}
\begin{Remark*} If $k_0$ is a product kernel of the form $k_0(z_1, \ldots, z_d)
= k_{00}(z_1) \cdots\break k_{00}(z_d)$, and
\[
k_{00}(x) \leq\gamma_2 e^{-\gamma_3 x^2 } \qquad\mbox{for every $x \in
R$,}
\]
then condition (K2) holds.


Assume the above conditions, then it is possible to control the
difference between $\hat{V}^*(z_k)$ and $V^*(z_k)$ using the following result.
\end{Remark*}
\begin{lm} \label{lm:normal}
Suppose that conditions \textup{(R1)--(R3)} and \textup{(K1)--(K2)} hold. Suppose that
$f_{n, 1}, \ldots, f_{n, k_n}$ are functions defined on $\mathcal
{X}\times
\mathcal{Y}\times\mathcal{Z}$, where $\mathcal{X}$, $\mathcal{Y}$
and $\mathcal{Z}$ are the ranges of $X$,
$Y$ and $Z$, respectively. Let $f_Z$ be the p.d.f.\vspace*{1pt} of $Z$, $\hat{f}_Z(z) =
(nh_n^{d})^{-1} \sum_{i=1}^{n} k_0(h_n^{-1}(z-Z_i))$ for $z\in
\mathcal{Z}$ and
$c_K = 1/\int k_0^2(s) \,ds$. For $z \in\mathcal{Z}$, let $w_i(z) =
n^{-1}h_n^{-d} w_{0,i}(z)/\hat{f}_Z(z)$ for $1 \leq i \leq n$ and
\[
W_{n, j} (z) = \sqrt{n h_n^d c_K f_Z(z) }  \Biggl(  \Biggl( \sum_{i=1}^{n}
w_i(z) f_{n, j } (X_i, Y_i,z)  \Biggr) - E\bigl(f_{n,j } (X, Y,z)|Z=z\bigr)
\Biggr)
\]
for $1 \leq j \leq k_n$.
Suppose that $\{ h_n \}_{n=1}^{\infty}$ and $\{ \varepsilon_n \}
_{n=1}^{\infty}$ are sequences of positive numbers such that
\[
c_{3,1} n^{-\alpha} \leq h_n \leq c_{3,2} n^{-\alpha}
\]
for some positive constants $c_{3,1}$ and $c_{3,2}$ and $1/(d+4) <
\alpha< 1/d$, and $h_n/\varepsilon_n = O(n^{-\beta})$ for some
$\beta
>0$. Let
%
%
\begin{equation} \label{eq:neighbor}
\mathcal{Z}(\varepsilon_n) =  \bigl\{ z \in\mathcal{Z}\dvtx \{ z' \in
R^d\dvtx \|z' - z \| <
\varepsilon_n \} \subset\mathcal{Z} \bigr\}
\end{equation}
and suppose that $z_1, \ldots, z_{n_Z}$ are points in $\mathcal{Z}
(\varepsilon_n)$ such that
%
%
\begin{equation} \label{eq:dense}
\| z_k - z_{k^*} \| \geq h_n \qquad\mbox{for $1 \leq k, k^* \leq n_Z$ and
$k\neq k^*$}
\end{equation}
for large $n$ and
%
%
\begin{equation} \label{eq:bdd}
{\max_{1\leq k \leq n_Z} \sup_{(x,y) \in\mathcal{X}\times\mathcal
{Y}}} |f_{n,j}(x,
y,z_k)| \leq C_n \qquad\mbox{for some $C_n \geq1$.}
\end{equation}
Suppose that $k_nn_ZC_n =O( (\ln n)^{1/16})$. Then there exist $W_{n,
1, j, k}$ and $W_{n, 2, j, k}\dvtx 1\leq j \leq k_n$, $1 \leq k \leq n_Z$
such that the joint distribution of $W_{n, 1, j,k} + W_{n, 2,j,k}$'s is
the same as the joint distribution of $W_{n, j}(z_k)$'s,
$\sum_{j=1}^{k_n} \sum_{k=1}^{n_Z} W^2_{n, 2,j,k} = O_P(\exp(-(\ln
n)^{1/9}))$, and $W_{n, 1, j,k}$'s are jointly normal with $EW_{n, 1,
j,k}=0$ and for $1 \leq j, \ell\leq k_n$ and $1 \leq k , k^* \leq n_Z$
\begin{eqnarray*}
&&\operatorname{Cov}(W_{n, 1, j, k}, W_{n, 1, \ell, k^*}) \\
&&\qquad = \cases{
\operatorname{Cov}\bigl(f_{n,j } (X, Y,z_k), f_{n, \ell} (X, Y,z_k)|Z=z_k\bigr), &\quad if $k
= k^*$; \cr
0, &\quad otherwise.}
\end{eqnarray*}
\end{lm}

The proof of Lemma \ref{lm:normal} is given in Section
\ref{sec:proof_lm_normal}.

The differences between $\hat{V}^*(z_k)$'s and $V^*(z_k)$'s can be
controlled by applying Lemma \ref{lm:normal} and taking the
$f_{n,j}(X,Y,z)$'s to be the functions $\phi^*_{\ell}(X) \phi
^*_{\ell
'}(X)$, $\phi^*_{\ell}(X) \psi^*_m(Y)$ and $\psi^*_m(Y) \psi
^*_{m'}(Y)$, where $1 \leq\ell\leq\ell' \leq p_n$ and $1 \leq m
\leq
m' \leq q_n$. In such case, (\ref{eq:bdd}) holds under the following
conditions.
\begin{enumerate}[(B2)]
\item[(B1)] For each $(p,q)$, $| \phi_{p, k} | \leq1$ and $|\psi_{q,
\ell} | \leq1$ for $1 \leq k \leq p$ and $1 \leq\ell\leq q$.
\item[(B2)] There exists $\{ \delta_n \}$: a sequence of positive
numbers such that for $1 \leq k \leq n_Z$, the smallest eigenvalues
of the matrices $V_{1,1}(z_k)$ and $V_{2,2}(z_k)$ are greater than or
equal to $\delta_n$.
\end{enumerate}
Under the above conditions, the $\hat{\rho}(z_k)$'s are consistent, as
stated in Theorem \ref{thm:consistency}.
\begin{theorem} \label{thm:consistency}
Suppose that (\ref{eq:app1}), (\ref{eq:app2}), (\ref{eq:basis}),
conditions \textup{(R1)--(R3)}, \textup{(K1)--(K2)} and \textup{(B1)--(B2)}
hold. Suppose that $\{
h_n \}_{n=1}^{\infty}$ and $\{ \varepsilon_n \}_{n=1}^{\infty}$ are
sequences of positive numbers such that
\[
c_{3,1} n^{-\alpha} \leq h_n \leq c_{3,2} n^{-\alpha}
\]
for some positive constants $c_{3,1}$ and $c_{3,2}$ and $1/(d+4) <
\alpha< 1/d$, and $h_n/\varepsilon_n = O(n^{-\beta})$ for some
$\beta
>0$. Suppose that $z_1, \ldots, z_{n_Z}$ are points in $\mathcal{Z}
(\varepsilon_n)$ [defined in (\ref{eq:neighbor})] such that (\ref
{eq:dense}) holds and
%
%
\begin{equation} \label{eq:small_nzpq}
n_Z (p_n+q_n)^2 \max\{ 1, \delta_n^{-1}(p_n+q_n) \} = O( (\ln n)^{1/16}).
\end{equation}
Then
%
%
\begin{equation} \label{eq:consistency1}
\sum_{k=1}^{n_Z}  \bigl( \hat{\rho}^2(z_k) - \rho_{p_n, q_n}^2(z_k)
 \bigr)^2 = O_P((nh_n^d)^{-1} (\ln n)^{1/4})
\end{equation}
and
%
%
\begin{equation} \label{eq:consistency2}\qquad
\Biggl( \sum_{k=1}^{n_Z} \hat{f}_Z(z_k) \hat{\rho}^2 (z_k) - \sum
_{k=1}^{n_Z} f_Z(z_k) \rho_{p_n, q_n}^2 (z_k)  \Biggr)^2
= O_P \biggl( \frac{ (\ln n)^{5/16} }{nh_n^d}  \biggr).
\end{equation}
\end{theorem}

The proof of Theorem \ref{thm:consistency} is given in Section
\ref{sec:proof_thm_consistency}.

The next result deals with the asymptotic distribution of $\sum
_{k=1}^{n_Z} \hat{f}_Z(z_k) \hat{\rho}^2 (z_k)$ when $X$ and $Y$ are
conditionally independent given $Z$.
\begin{theorem} \label{thm:distribution}
Suppose that the conditions in Theorem \ref{thm:consistency} hold and
$X$ and $Y$ are conditionally independent given $Z$. Then there exist
random variables $\tilde{f}_k$, $\tilde{\rho}^2(z_k)$ and
$\lambda_k\dvtx
1 \leq k \leq n_Z$ such that $\sum_{k=1}^{n_Z} \tilde{f}_k \tilde
{\rho
}^2(z_k)$ has the same distribution as $\sum_{k=1}^{n_Z} \hat{f}_Z(z_k)
\hat{\rho}^2 (z_k)$ and
\[
nh_n^d c_K \sum_{k=1}^{n_Z} \tilde{f}_k \tilde{\rho}^2(z_k) - \sum
_{k=1}^{n_Z} \lambda_k = O_P(\exp(-0.5 (\ln n)^{1/9}) (\ln n)^{3/32} ),
\]
where the $\lambda_k$'s are independent and each $\lambda_k$ has the
same distribution as the largest eigenvalue of a matrix $CC^T$, where
$C$ is a $(p_n -1) \times(q_n-1)$ matrix whose elements are i.i.d. $N(0,1)$.
\end{theorem}

The proof of Theorem \ref{thm:distribution} is given in Section
\ref{sec:proof_thm_distribution}. The result in Theorem \ref
{thm:distribution} is similar to that in Lemma 7.2 in
\cite{MR1647653}. 
The difference is that the asymptotic result here is derived as the
sample size $n$, $p_n$ and $q_n$ all tend to $\infty$, while in \cite
{MR1647653}, 
the result is derived as $n$ tends to $\infty$, but $p_n$ and $q_n$ are
held fixed.

Theorem \ref{thm:distribution} suggests the test that rejects the
conditional independence hypothesis at approximate level $a$ if
%
%
\begin{equation} \label{eq:test1}
nh_n^d c_K \sum_{k=1}^{n_Z} \hat{f}_Z(z_k) \hat{\rho}^2 (z_k) >
F_{n_Z, p, q}^{-1}(1-a),
\end{equation}
where $F_{n_Z, p, q}$ is the cumulative distribution function of $\sum
_{k=1}^{n_Z} \lambda_k$.

One can estimate $F_{n_Z, p, q}^{-1}(1-a)$ in (\ref{eq:test1}) using
simulated data, but it is also possible to use a normal approximation.
Since the $\lambda_k$'s are i.i.d., the central limit theorem suggests
the asymptotic normality of $\sum_{k=1}^{n_Z} \lambda_k$ and $\sum
_{k=1}^{n_Z} \hat{f}_Z(z_k) \hat{\rho}^2 (z_k)$. The following
corollary gives the conditions that guarantee the asymptotic normality
of $\sum_{k=1}^{n_Z} \hat{f}_Z(z_k) \hat{\rho}^2 (z_k)$.
\begin{cor} \label{cor:normal}
Suppose that the conditions in Theorem \ref{thm:consistency} hold
%
%
\begin{equation} \label{eq:nz} \lim_{n\rightarrow\infty} \frac
{p_n^3 q_n^3 }
{ \sqrt{n_Z} (\max(p_n, q_n))^{1/3} } = 0
\end{equation}
and \textup{(i)} or \textup{(ii)} holds:
\begin{longlist}
\item $q_n = h(p_n)$, where $h$ is an increasing function such
that $\lim_{p\rightarrow\infty} h(p)/p$ exists and is greater than
or equal to 1.
\item $p_n = h(q_n)$, where $h$ is an increasing function such
that $\lim_{q\rightarrow\infty} h(q)/q$ exists and is greater than
or equal to 1.
\end{longlist}
Let $\mu_{p_n, q_n}$ and $\sigma^2_{p_n, q_n}$ be the mean and variance
of the largest eigenvalue of the matrix $CC^T$ in Theorem
\ref{thm:distribution}, respectively, and let the $\lambda_k$'s be as in
Theorem \ref{thm:distribution}, then
%
%
\begin{equation} \label{eq:cor_normal_var} \frac{ (\max(p_n,
q_n))^{1/6} }{ \sigma_{p_n, q_n} } = O(1)
\end{equation}
and
%
%
\begin{equation} \label{eq:lambda_normal}
\frac{ \sum_{k=1}^{n_Z} \lambda_k - n_Z \mu_{p_n,q_n} }{ \sqrt{ n_Z
\sigma^2_{p_n, q_n}} } \stackrel{\mathcal{D}}{\rightarrow} N(0,1)
\qquad\mbox{as $n \rightarrow
\infty$.}
\end{equation}
If $X$ and $Y$ are conditionally independent given $Z$, then
%
%
\begin{equation} \label{eq:stat_normal}\qquad
\frac{ nh_n^d c_K \sum_{k=1}^{n_Z} \hat{f}_Z(z_k) \hat{\rho}^2
(z_k) -
n_Z \mu_{p_n,q_n} }{ \sqrt{ n_Z \sigma^2_{p_n, q_n}} } \stackrel
{\mathcal{D}
}{\rightarrow} N(0,1) \qquad\mbox{as $n \rightarrow\infty$. }
\end{equation}
\end{cor}

The proof of Corollary \ref{cor:normal} is given in Section
\ref{sec:proof_cor_normal}. Corollary \ref{cor:normal} gives the test that
rejects the conditional independence hypothesis if
%
%
\begin{equation} \label{eq:reject}
\frac{ nh_n^d c_K \sum_{k=1}^{n_Z} \hat{f}_Z(z_k) \hat{\rho}^2
(z_k) -
n_Z \mu_{p_n,q_n} }{ \sqrt{ n_Z \sigma^2_{p_n, q_n}} } \geq\Phi^{-1}(1-a),
\end{equation}
where $\Phi$ is the cumulative distribution function for the standard
normal distribution. Here, $\mu_{p_n,q_n}$ and $\sigma^2_{p_n, q_n}$
can be approximated by the sample mean and variance of a random sample
from the distribution of the largest eigenvalue of the matrix $CC^T$.

To distinguish the two tests mentioned above, we will refer to the test
with rejection region in (\ref{eq:reject}) as test 1N and the test with
rejection region in (\ref{eq:test1}) as test~1. Note that under the
conditions in Corollary \ref{cor:normal}, test 1 does not differ from
test~1N much since the rejection region for test 1 can be written as
\[
\frac{ nh_n^d c_K \sum_{k=1}^{n_Z} \hat{f}_Z(z_k) \hat{\rho}^2
(z_k) -
n_Z \mu_{p_n,q_n} }{ \sqrt{ n_Z \sigma^2_{p_n, q_n}} } \geq I + \Phi
^{-1}(1-a),
\]
where
%
%
\begin{equation} \label{eq:test_diff}
I = \frac{F_{n_Z, p, q}^{-1}(1-a) - n_Z \mu_{p_n, q_n}}{\sqrt{ n_Z
\sigma^2_{p_n, q_n}}} - \Phi^{-1}(1-a) = o(1)
\end{equation}
by (\ref{eq:lambda_normal}).
Therefore, both tests 1 and 1N are of asymptotic significance level
$a$. Below we will discuss the consistency and asymptotic power of test
1N only since the same properties of test 1 can be established
similarly using (\ref{eq:test_diff}).

Suppose all the conditions in Theorem \ref{thm:consistency} hold, then
test that 1N is also consistent if the $z_k$'s are chosen in a way such that
there exist a constant $c_3 >0$ and a sequence $\{ \eta_{1,n} \}
_{n=1}^{\infty}$ such that $\eta_{1,n}>0$ for every $n$, $\lim
_{n\rightarrow
\infty} \eta_{1,n} =0$ and
%
%
\begin{equation} \label{eq:uniform}
\frac{1}{n_Z} \sum_{k=1}^{n_Z} f_Z(z_k) \rho^2_{p_n,q_n} (z_k) - c_3
E\rho_{p_n,q_n}^2(Z) = o_P(\eta_{1,n}).
\end{equation}
To see that test 1N is consistent, note that $0 \leq\mu_{p_n, q_n}
\leq E\operatorname{tr}(CC^T)$ and $\sigma^2_{p_n, q_n} \leq E(\operatorname{tr}(CC^T))^2$,
where $CC^T$ is as in Theorem \ref{thm:distribution}. Therefore, $\mu
_{p_n, q_n} = O(p_nq_n)$ and $\sigma^2_{p_n, q_n}
= O(p_n^2 q_n^2)$. Then it follows from (\ref{eq:consistency2}), (\ref
{eq:uniform}) and Fact \ref{fact:app} that
$n_Z^{-1} \sum_{k=1}^{n_Z} \hat{f}_Z(z_k) \hat{\rho}^2 (z_k) - c_3
E\rho
_1^2(X,Y|Z)= O_P((\ln n)^{5/32}/n_Z\sqrt{nh_n^d})
+ o_P(\eta_{1,n}) + c_3 E\rho_{p_n,q_n}^2(Z) - c_3 E\rho_1^2(X,Y|Z)
= o_P(1)$, so
\begin{eqnarray*}
&& \frac{ nh_n^d c_K \sum_{k=1}^{n_Z} \hat{f}_Z(z_k) \hat{\rho}^2 (z_k)
- n_Z \mu_{p_n,q_n} }{ \sqrt{ n_Z \sigma^2_{p_n, q_n}} } \\
&&\qquad \geq \frac{ \sqrt{n_Z}  ( nh_n^d c_K (c_3 E\rho_1^2(X,Y|Z) +
o_P(1)) + O(p_nq_n)  )}{ c_{2,1} p_n q_n},
\end{eqnarray*}
where $c_{2,1} >0$ is a constant. Thus, the left-hand side in (\ref
{eq:reject}) tends to $\infty$ as $n \rightarrow\infty$ when $E\rho
_1^2(X,Y|Z) >0$, which implies that the probability that (\ref
{eq:reject}) holds tends to 1 if $X$ and $Y$ are not conditionally
independent given $Z$.\vspace*{1pt}

Test 1N can also reject an alternative where $E\rho_{p_n,q_n}^2(Z)$ is
small under the conditions in Theorem \ref{thm:consistency}. Indeed,
for $\{ \eta_{1,n} \}_{n=1}^{\infty}$ such that $\eta_{1,n}>0$ for
every $n$, $\lim_{n\rightarrow\infty} \eta_{1,n} =0$ and (\ref{eq:uniform})
holds, if
%
%
\begin{equation} \label{eq:alternative1}
\frac{\max ( \eta_{1,n}, { (\ln n)^{5/32} }/({n_Z\sqrt
{nh_n^d}})) }{E\rho_{p_n,q_n}^2(Z)} =o(1),
\end{equation}
then the probability that (\ref{eq:reject}) holds tends to 1 since
\begin{eqnarray*}
&& \frac{ nh_n^d c_K \sum_{k=1}^{n_Z} \hat{f}_Z(z_k) \hat{\rho}^2 (z_k)
- n_Z \mu_{p_n,q_n} }{ \sqrt{ n_Z \sigma^2_{p_n, q_n}} } \\
&&\qquad \geq \biggl( \sqrt{n_Z}  \biggl( nh_n^d c_K  \biggl(c_3 E\rho
_{p_n,q_n}^2(Z)\\
 &&\hspace*{102.4pt}{} + O_P \biggl( \frac{ (\ln n)^{5/32} }{n_Z\sqrt{nh_n^d}}
 \biggr) + o_P(\eta_{1,n})  \biggr)
+ O(p_nq_n)  \biggr)\biggr)\\
&&\qquad\quad{}\times({ c_{2,1} p_n q_n})^{-1},
\end{eqnarray*}
where $p_nq_n/(nh_n^dE\rho_{p_n,q_n}^2(Z)) =O((\ln
n)^{1/16}/(n_Znh_n^dE\rho_{p_n,q_n}^2(Z))) = o(1)$ by (\ref
{eq:small_nzpq}) and (\ref{eq:alternative1}), and $p_nq_n/(\sqrt{n_Z}
nh_n^dE\rho_{p_n,q_n}^2(Z)) = o(1)$. In summary, test 1N can reject an
alternative where $E\rho_{p_n,q_n}^2(Z)$ tends to zero at a rate that
is slower than $\max( \eta_{1,n}, (\ln n)^{5/32}/(n_Z\sqrt{nh_n^d}))$,
where $\eta_{1,n}$ is determined by (\ref{eq:uniform}). An example that
satisfies (\ref{eq:uniform}) and the conditions in Corollary \ref
{cor:normal} will be given in Section~\ref{sec:ex}. In that example,
$\eta_{1,n} = p_n^{11}n_Z^{-1/d}$.

\section{An example} \label{sec:ex}
In this section, an example is given to illustrate the verification of
the conditions in Corollary \ref{cor:normal}, assuming (R1)--(R3) and
the condition that there exists a positive constant $c_{1,1}$ such that
%
%
\begin{eqnarray} \label{eq:bdd_below}
f_{X|Z}(x|z) \geq c_{1,1} \quad\mbox{and}\quad f_{Y|Z}(y|z) \geq c_{1,1}
\nonumber\\[-8pt]\\[-8pt]
\eqntext{\mbox{for all } (x,y,z)
\in\mathcal{X}\times\mathcal{Y}\times\mathcal{Z},}
\end{eqnarray}
where
$f_{X|Z}(\cdot|z)$ and $ f_{Y|Z}(\cdot|z)$ are conditional probability
densities of $X$ and $Y$, respectively, given $Z=z$, with respect
to Lebesgue measures.
\begin{ex}\label{example1}
Suppose that $X$, $Y$ and $Z$ are random vectors that take values in
$[0,1]^{d_x}$, $[0,1]^{d_y}$ and $[0,1]^{d}$, respectively.
Suppose that \textup{(R1)--(R3)}, and (\ref{eq:bdd_below}) hold. Choose the
basis functions as follows. Let $\Lambda$ be the set of all
positive integers and $\Lambda(k) = \{ m^k\dvtx m \in\Lambda\}$ for $k
\in\Lambda$.
For $k$, $i_1, \ldots, i_k \in\Lambda$ and $h_0 >0$, let
\[
h_{k, h_0, i_1, \ldots, i_k}(x_1, \ldots, x_k) = \prod_{j=1}^k
I_{A_{i_j, h_0}}(x_j)\qquad
\mbox{for $(x_1, \ldots, x_k) \in[0,1]^k$,}
\]
where
\[
A_{i_j, h_0} =  \cases{
(h_0(i_j -1) , h_0 i_j], &\quad if $i_j > 1$;\cr
[h_0(i_j -1) , h_0 i_j], &\quad if $i_j =1$.}
\]
For $p$, $q$, $r \in\Lambda$, let
\begin{eqnarray*}
\{ \phi_{p,i}\dvtx 1 \leq i \leq p \} &=& \{ h_{d_x, p^{-1/d_x} , i_1,
\ldots, i_{d_x} }\dvtx 1 \leq i_1, \ldots, i_{d_x} \leq p^{1/d_x} \},
\\
\{ \psi_{q,j}\dvtx 1 \leq j \leq q \} &=& \{ h_{d_y, q^{-1/d_y} , i_1,
\ldots
, i_{d_y} }\dvtx 1
\leq i_1, \ldots, i_{d_y} \leq q^{1/d_y} \}
\end{eqnarray*}
and
\[
\{ \theta_{r, k}\dvtx 1 \leq k \leq r \} = \{ h_{d, r^{-1/d} , i_1, \ldots,
i_d }\dvtx 1 \leq
i_1, \ldots, i_d \leq r^{1/d} \}.
\]
Take $k_0$ to be the product kernel function such that
\[
k_0(z_1, \ldots, z_d) = k_{00}(z_1) \cdots k_{00}(z_d),
\]
where $k_{00}$ is the probability density function for
the standard normal distribution. Let $h_n = n^{-a}$, where $1/(d+4) <
a < 1/d$. Let $n_Z^*$ to be the largest number
in $\Lambda(d)$ such that $n_Z^* \leq(\ln n)^{1/32}$, and let
\[
\{ z_k\dvtx 1 \leq k \leq n_Z \} =  \biggl\{  \biggl( \frac{i_1}{(n_Z^*)^{1/d}},
\ldots, \frac{i_d}{(n_Z^*)^{1/d}}  \biggr)\dvtx 1 \leq i_1, \ldots, i_d <
(n_Z^*)^{1/d}  \biggr\},
\]
so $n_Z = ((n_Z^*)^{1/d}-1)^d$.
Suppose that $\{ p_n \}$ is a sequence in $\Lambda(d_x) \cap\Lambda
(d_y)$ such that $\lim_{n\rightarrow\infty} p_n = \infty$ and $q_n
= p_n$. If
%
%
\begin{equation} \label{eq:pn1}
p_n^{12} \leq n_Z,
\end{equation}
then all the conditions in Corollary \ref{cor:normal} hold. If
%
%
\begin{equation} \label{eq:pn2}
p_n^{12} \leq n_Z^{1/d},
\end{equation}
then (\ref{eq:uniform}) holds with $\eta_{1,n} = p_n^{11}n_Z^{-1/d}$.
\end{ex}
\begin{pf}
We will first show that all the conditions in Corollary \ref
{cor:normal} hold assuming (\ref{eq:pn1}).
It is clear that (\ref{eq:app1}), (\ref{eq:app2}) and (\ref
{eq:basis}), and conditions (B1), (K1) and (K2) hold.

To find the $\delta_n$ in condition (B2), note that for $z \in
\mathcal{Z}$,
the smallest eigenvalue of
$V_{1,1}(z)$ is the minimum of $\{ E(\phi_{p_n, i}(X)|Z=z)\dvtx 1 \leq i
\leq p_n \}$, which is the minimum of
$\{ E(h_{ d_x, p_n^{-1/{d_x}} , i_1, \ldots, i_{d_x} } (X)|Z=z)\dvtx 1
\leq
i_1, \ldots, i_{d_x} \leq p_n^{1/{d_x}} \}$. Under (\ref
{eq:bdd_below}), for $m \in\Lambda$ and $1 \leq i_1, \ldots, i_{d_x}
\leq m$,
\begin{eqnarray*}
&& E\bigl(h_{d_x, 1/m , i_1, \ldots, i_{d_x}} (X)|Z=z\bigr) \\
&&\qquad = \int_{(i_1-1)/m}^{i_1/m} \cdots\int_{(i_{d_x}-1)/m}^{i_{d_x}/m}
f_{X|Z} (x_1, \ldots, x_{d_x}|z) \,dx_{d_x} \cdots d x_1
\geq\frac{ c_{1,1} }{ m^{d_x}}.
\end{eqnarray*}
Take $m=p_n^{1/{d_x}}$, and we have that the smallest eigenvalue of $V_{1,1}(z)$
is at least $c_{1,1}/p_n$. Similarly, $c_{1,1}/p_n$ is also a lower
bound for the smallest eigenvalue of $V_{2,2}(z)$ and (B2) holds
with $\delta_n = c_{1,1}/p_n$. Furthermore, (\ref{eq:small_nzpq})
holds since
\[
n_Z (p_n+q_n)^2 \max\{ 1, \delta_n^{-1}(p_n+q_n) \} = O( n_Z p_n^4) =
O(n_Z^2).
\]

Finally, the $z_k$'s are in $\mathcal{Z}(\varepsilon_n)$ with
$\varepsilon_n =
(n_Z^*)^{-1/d}$ and $h_n/\varepsilon_n = O(n^{-\beta})$ for $0 <
\beta
< \alpha$.
For $1 \leq k, k^* \leq n_Z$, and $k\neq k^*$, $\| z_k - z_{k^*} \|
\geq(n_Z^*)^{-1/d} \geq n^{-a}$, so (\ref{eq:dense}) holds. Also,
(\ref
{eq:nz}) holds since
\[
\frac{p_n^3 q_n^3 } { \sqrt{n_Z} (\max(p_n, q_n))^{1/3} } = p_n^{-1/3}
\sqrt{ \frac{p_n^{12}}{n_Z} } =o(1).
\]
Therefore, all the conditions in Corollary \ref{cor:normal} hold for
this example.


The verification of (\ref{eq:uniform}) is based on the fact that there
exist positive constants $c_{4,1}$ and $\eta_0$ such that
%
%
\begin{equation} \label{eq:lipschitz}
| \rho^2_{p_n,q_n}(z) - \rho^2_{p_n,q_n}(z') | \leq c_{4,1} p_n ^{11}
\| z - z' \| \qquad\mbox{if $p_n^3 \| z - z' \| < \eta_0$.}
\end{equation}
Below we will first check (\ref{eq:uniform}) assuming that
(\ref{eq:lipschitz}) holds and then prove (\ref{eq:lipschitz}).
Suppose that (\ref{eq:pn2}) holds. Let $g_n(z) = f_Z(z)\rho^2_{p_n,
q_n}(z)$. Since $f_Z$ is Lipschitz continuous, (\ref{eq:lipschitz})
implies that there exists a constant $c_{4,2} > 0$ such that
\[
| g_n(z) - g_n(z') | \leq c_{4,2} p_n^{11} \| z - z' \|\qquad
\mbox{if $p_n^3 \| z - z' \| < \eta_0$.}
\]
Let $\{ z_{1+n_Z}, \ldots, z_{n_Z^*} \}$ be the set
\[
\biggl\{  \biggl( \frac{i_1}{(n_Z^*)^{1/d}}, \ldots, \frac
{i_d}{(n_Z^*)^{1/d}}  \biggr)\dvtx
1 \leq i_1, \ldots, i_d \leq(n_Z^*)^{1/d}  \biggr\} \cap\{ z_k\dvtx 1
\leq k \leq n_Z \}^c,
\]
then
\[
\Biggl| \sum_{k=1}^{n_Z^*} g_n(z_k)  \biggl( \frac{1}{(n_Z^*)^{1/d}}
 \biggr)^d
- \int_{\mathcal{Z}} g_n(z) \,dz \Biggr| \leq2 c_{4,2} p_n^{11} \sqrt
{d}\biggl( \frac{1}{n_Z^*}  \biggr)^{1/d},
\]
if $p_n^3 (n_Z^*)^{-1/d} < \eta_0$.
Since $|g_n(z)| \leq c_0$ by (R3) and there exists a positive constant
$c_{4,3}$ depending on $d$ such that
\[
n_Z^* - n_Z
\cases{\leq c_{4,3} (n_Z^*)^{1/d}, &\quad if $d \geq2$; \cr
=1, &\quad if $d =1$,}
\]
we have
\begin{eqnarray*}
&&  \Biggl| n_Z^{-1} \sum_{k=1}^{n_Z} f_Z(z_k) \rho^2_{p_n,q_n} (z_k) -
\frac{\int_{\mathcal{Z}} f_Z(z) \rho^2_{p_n,q_n} (z) \,dz}{ \int
_{\mathcal{Z}} 1 \,dz}
 \Biggr| \\
&&\qquad =  \Biggl| \frac{n_Z^*}{n_Z}  \Biggl( \frac{1}{n_Z^*} \sum
_{k=1}^{n_Z^*} g_n(z_k) - \int_{\mathcal{Z}} g_n(z) \,dz  \Biggr)
\\
&&\qquad\quad\hspace*{1.6pt}{} - \frac{\sum_{k=1+n_Z}^{n_Z^*} g_n(z_k) }{n_Z} + \biggl( \frac
{n_Z^*}{n_Z} - 1  \biggr) \int_{\mathcal{Z}} g_n(z) \,dz  \Biggr| \\
&&\qquad \leq \frac{n_Z^*}{n_Z}  \Biggl| \frac{1}{n_Z^*} \sum_{k=1}^{n_Z^*}
g_n(z_k) - \int_{\mathcal{Z}} g_n(z) \,dz  \Biggr| + c_0  \biggl( 1 +
\int_{\mathcal{Z}} 1\,dz
 \biggr)  \biggl(
\frac{n_Z^* - n_Z}{n_Z}  \biggr) \\
&&\qquad\leq\frac{ c_{4,4} p_n^{11} }
{n_Z^{1/d} }
\end{eqnarray*}
for some constant $c_{4,4} > 0$ if $p_n^3 (n_Z^*)^{-1/d} < \eta_0$.
Since $p_n^{12} \leq n_Z^{1/d}$, $p_n^{3}n_Z^{-1/d} =o(1)$, so
\[
\Biggl| n_Z^{-1} \sum_{k=1}^{n_Z} f_Z(z_k) \rho^2_{p_n,q_n} (z_k) -
\frac
{\int_{\mathcal{Z}} f_Z(z) \rho^2_{p_n,q_n} (z) \,dz}{ \int_{\mathcal
{Z}} 1 \,dz}  \Biggr| =
O_P \biggl( \frac{ p_n^{11} }{ n_Z^{1/d} }  \biggr)
\]
and $p_n^{11}n_Z^{-1/d} =o(1)$. Take $\eta_{1,n} = p_n^{11}n_Z^{-1/d}$
and $c_3= (\int_{\mathcal{Z}} 1 \,dz)^{-1} = 1$,
then (\ref{eq:uniform}) holds.

It remains to prove (\ref{eq:lipschitz}). Recall that for $z\in
\mathcal{Z}$,
$\rho^2_{p_n,q_n}(z)$ is the largest eigenvalue of $g(V(z), a_{n,*})$,
as mentioned in Section
\ref{section:estimator}. Thus, $| \rho^2_{p_n,q_n}(z) - \rho
^2_{p_n,q_n}(z')|$ is bounded by $\| g(V(z), a_{n,*}) - g(V(z'),
a_{n,*}) \|$.
For $1 \leq i, j \leq2$, let $g_{i,j}^*$ be as defined in (\ref
{eq:gstar}) and let
$\Delta_{i,j}=g_{i,j}^*(V(z')) -g_{i,j}^*(V(z))$ for $1 \leq i, j \leq
2$, then from the fact that $\| AB \| \leq\| A \| \| B \|$ for two
matrices $A$ and $B$, we have
%
%
\begin{eqnarray}\label{eq:consistent3m}
&& \| g(V(z), a_{n,*}) - g(V(z'), a_{n,*}) \| \nonumber\\
&&\qquad \leq \prod_{i=1}^2 \prod_{j=1}^2  \bigl( \| g^*_{i,j}(V(z)) \| +
\|
\Delta_{i,j} \|  \bigr) - \prod_{i=1}^2 \prod_{j=1}^2 \|
g^*_{i,j}(V(z)) \| \\
&&\qquad\quad{} + \| g_{1,1}(V(z')) -g_{1,1}(V(z)) \| \| a_{n,*} \|^2.\nonumber
\end{eqnarray}
The bounds for the $\| g^*_{i,j}(V(z)) \|$'s are derived as follows.
Since the elements in $V(z)$ are bounded by 1
and the smallest eigenvalue of $g_{i,i}(V(z))$ is at least
$c_{1,1}/p_n$ for $1 \leq i \leq2$, we have
\begin{eqnarray*}
\max( \| g^*_{1,2}(V(z)) \|, \| g_{2,1}^*(V(z)) \| ) &\leq& p_n,
\\
\| g_{1,1}^*(V(z)) \|^2 &\leq& \frac{p_n^2}{ ( c_{1,1}/p_n)^2} = \frac
{p_n^4}{c_{1,1}^2}
\end{eqnarray*}
and
\[
\| g_{2,2}^*(V(z)) \| \leq\frac{p_n^2}{c_{1,1}}.
\]
To find bounds for $\| g_{1,1}(V(z')) -g_{1,1}(V(z)) \| $ and $\|
\Delta
_{i,j} \|$'s, note that from (R3), each element in $g_{i,j}(V(z'))
-g_{i,j}(V(z))$ is bounded by $\sqrt{d} \int h(x,y) \,d\mu(x,y)\| z-
z' \|$, so
\begin{eqnarray*}
&& \max\bigl( \| \Delta_{1,2} \|, \|\Delta_{2,1} \|, \| g_{1,1}(V(z'))
-g_{1,1}(V(z)) \| \bigr) \\
&&\qquad \leq p_n \sqrt{d} \int h(x,y) \,d\mu(x,y) \| z- z' \|.
\end{eqnarray*}
For $1 \leq i \leq2$, by Fact \ref{fa:matrix},
\[
\| \Delta_{i,i} \| \leq\frac{\| g_{i,i}^*(V(z)) \|^2 \|g_{i,i}(V(z'))
- g_{i,i} (V(z)) \|}{1-\| g_{i,i}^*(V(z)) \| \|g_{i,i}(V(z')) - g_{i,i}
(V(z)) \|},
\]
if $\| g_{i,i}^*(V(z)) \| \|g_{i,i}(V(z')) - g_{i,i} (V(z)) \| < 1$, so
\[
\| \Delta_{i,i} \| \leq\frac{2\sqrt{d}p_n^5}{c_{1,1}^2} \int h(x,y)\,
d\mu(x,y) \| z- z' \|,
\]
if
%
%
\begin{equation} \label{eq:smalldelta} \frac{\sqrt{d}p_n^3}{c_{1,1}}
\int h(x,y) \,d\mu(x,y) \| z- z' \| < \frac{1}{2}.
\end{equation}
To give a bound for $\| a_{n,*} \|$, note that the smallest eigenvalue
of $g_{1,1}(V(z))$ is at least $c_{1,1}/p_n$ and at most
\[
\frac{a_{n,*}^T g_{1,1}(V(z)) a_{n,*}}{a_{n,*}^T a_{n,*}} = \frac
{1}{\|
a_{n,*} \|^2},
\]
so
\[
\| a_{n,*} \| \leq\sqrt{\frac{p_n}{c_{1,1}}}.
\]
From (\ref{eq:consistent3m}) and the above bounds for $\| a_{n,*} \|$,
the $\| g^*_{i,j}(V(z)) \|$'s and $\| \Delta_{i,j} \|$'s, we have
\[
\| g(V(z), a_{n,*}) - g(V(z'), a_{n,*}) \| \leq c_{4,1} p_n^{11} \|
z-z' \|
\]
for some constant $c_{4,1}$ if (\ref{eq:smalldelta}) holds. Therefore,
(\ref{eq:lipschitz})
holds and the proof for the results in Example \ref{example1} is complete.
\end{pf}

\section{Simulation studies} \label{sec:sim}
In this section, results of several simulation experiments are
presented. Those experiments are designed to demonstrate the
performance of test 1 introduced in Section \ref{sec:test}.

In Section \ref{sec:test}, test 1N is also introduced, but no
simulation studies are done for it in this section. The reason is as
follows. Test 1N is constructed based on the normal approximation for
$\sum_{k=1}^{n_Z} \lambda_k$. Using the parameter set-up in Table
\ref{ta:parameter}, the selected $n_Z$ is only 4 or 5 and the normal
approximation for $\sum_{k=1}^{n_Z} \lambda_k$ is not expected to
work well.

For simplicity, in all the simulation experiments here, $X$, $Y$, $Z$
are one dimensional and only the following distributions for $(X,Y,Z)$
are considered.
\begin{enumerate}[(M3)]
\item[(M1)] $(X, Y) = (\Phi(Z \epsilon_1), \Phi(Z \epsilon_2))$, where
$\epsilon_1$, $\epsilon_2$ and $Z$ are independent, $Z$ follows the
uniform distribution on $[0,1]$, and $\epsilon_i$ follows the standard
normal distribution for $i=1$, 2.

\item[(M2)]
$Z$ follows the standard normal distribution, and the conditional
distribution of $(X, Y)$ given $Z=z$ is bivariate normal with mean $\mu
$ and covariance matrix $\Sigma$, where
%
%
\begin{equation} \label{eq:mu_sigma}
\mu=  \pmatrix{
0 \cr 0},\qquad
\Sigma= \pmatrix{1 & \rho(z) \cr \rho(z) & 1}
\end{equation}
and the $\rho(z)$ in (\ref{eq:mu_sigma}) is taken to be $a(|1-2\Phi
(z)|)$ with $a \in\{ 0, 0.1, 0.3 \}$.

\item[(M3)] $(X,Y,Z) = (\Phi(X_0), \Phi(Y_0), \Phi(Z_0))$, where $Z_0$
follows the $t$-distribution with degree of freedom 1, and the
conditional distribution of $(X_0, Y_0)$ given $\Phi(Z_0)=z$ is
bivariate normal with mean $\mu$ and covariance matrix $\Sigma$, where
$\mu$ and $\Sigma$ are as in (\ref{eq:mu_sigma}) and the $\rho(z)$ in
(\ref{eq:mu_sigma}) is taken to be $a(|1-2z|)$ with $a \in\{ 0, 0.1,
0.3 \}$.
\end{enumerate}
Here, (M1) is used for parameter selection and (M2) and (M3) are used
for checking the power of test 1. In (M1), $X$ and $Y$ are
conditionally independent given $Z$.
In (M2) and (M3), $\rho_1(X,Y|Z=z) = \rho(z)$ and $E\rho_1(X,Y|Z)$ is
proportional to $a$.

The details of parameter selection are given in Section \ref{sec:para}
and the experimental results are given in Section \ref{sec:exp}.

\subsection{Parameter selection} \label{sec:para}
To apply test 1, certain parameters need to be chosen, including the
kernel function $k_0$, the kernel bandwidth $h_n$, the basis functions
$\phi_{p_n,i}$'s and $\psi_{q_n,j}$'s and the evaluation points
$z_k$'s, which are chosen as follows.
\begin{enumerate}[(S1)]
\item[(S1)] $k_0$ and the basis functions $\phi_{p,i}$'s and $\psi
_{q,j}$'s are chosen as in Example \ref{example1} in Section \ref{sec:ex} with
$p_n = q_n =2$. Since the basis functions are supported on $[0,1]$,
if $X$, $Y$ and $Z$ do not take values in $[0,1]$ [such as in (M2)],
then the data $\{ (X_i, Y_i, Z_i) \}_{i=1}^n$ will be transformed to
$\{ (\Phi(X_i), \Phi(Y_i), \Phi(Z_i)) \}_{i=1}^n$ before applying test
1. The bandwidth $h_n$ is chosen to be the $h$ that minimizes
%
%
\begin{equation} \label{eq:s1}
\int_{ 0.143h^{0.121} }^{1- 0.143h^{0.121}} E \bigl( \hat{f}_Z(z) - 1
\bigr)^2 \,dz
\end{equation}
over $(0, 0.5]$, where $\hat{f}_Z$ is the kernel density estimator
based on a sample of size $n$ from the uniform distribution on $[0,1]$
with kernel $k_0$ and bandwidth~$h$. Below are the $h_n$'s used for
different $n$'s.

The $z_k$'s are points in $I_n = [0.143h_n^{0.121}, 1 -
0.143h_n^{0.121}]$ such that $z_k = 0.143h_n^{0.121} + (k-1)h_{0,n}$,
where $h_{0,n}$ is a given positive number. Here, the $\varepsilon_n$
is taken to be $0.143h_n^{0.121}$, so the $z_k$'s are chosen so that
they are $0.143h_n^{0.121}$ away from the boundary and the integral in
(\ref{eq:s1}) is over $[0.143h^{0.121}, 1 - 0.143h^{0.121}]$.
\end{enumerate}
With the parameter set-up in (S1), it remains to choose $h_{0,n}$.
The $h_{0,n}$ is chosen to be the smallest multiple of 0.01 such that
the distribution for the test 1 statistic $nh_n^d c_K \sum_{k=1}^{n_Z}
\hat{f}_k \hat{\rho}^2(z_k)$ based on 1000 samples of size $n$ from
(M1) is similar to the distribution of $\sum_{k=1}^{n_Z} \lambda_k$
($\chi^2$ with $n_Z$ degrees of freedom), as stated in Theorem \ref
{thm:distribution}. The one-sample Kolmogorov--Smirnov test is used to
determine whether the two distributions are similar. Below are the
$h_{0,n}$'s used for $n=10\mbox{,}000$ and $n=5000$.

For the above procedure for selecting $h_{0,n}$, when $n=500$ or
$n=1000$, it seems that the distribution of $nh_n^d c_K \sum
_{k=1}^{n_Z} \hat{f}_k \hat{\rho}^2(z_k)$ cannot be approximated well
by the distribution of $\sum_{k=1}^{n_Z} \lambda_k$, regardless what
$h_{0,n}$ is used. To overcome this problem, one may use local
bootstrap to determine the rejection region.

The idea of using local bootstrap is to draw samples $\{ (X^*_i, Y^*_i,
Z^*_i) \}_{i=1}^n$ from the distribution of $(X^*, Y^*, Z^*)$, where
$Z^*$'s distribution is close to the distribution of $Z$ and the
conditional distributions of $X^*$ given $Z^*=z$ and $Y^*$ given
$Z^*=z$ are close to the conditional distributions of $X$ given $Z=z$
and $Y$ given $Z=z$, yet $X^*$ and $Y^*$ are conditionally independent
given $Z^*$. Therefore, if $X$ and $Y$ are conditionally independent
given $Z$, then the local bootstrap resamples $\{ (X^*_i, Y^*_i, Z^*_i)
\}_{i=1}^n$ should behave like a random sample from $(X,Y,Z)$. One can
then compute the test~1 statistic $nh_n^d c_K \sum_{k=1}^{n_Z} \hat
{f}_k \hat{\rho}^2(z_k)$ for the original sample and for each local
bootstrap resample. If the statistic computed based on the original
sample is larger than $(1-a)$\% of the statistics computed based on the
local bootstrap resamples, then the conditional independence hypothesis
is rejected at level~$a$.

%
%
\begin{table}
\tablewidth=230pt
\caption{Selected $h_n$'s for different $n$'s} \label{ta:hn}
\begin{tabular*}{\tablewidth}{@{\extracolsep{\fill}}lcccc@{}}
\hline
$\bolds n$ & \textbf{10\mbox{,}000} & \textbf{5000} & \textbf{1000} & \textbf{500} \\
\hline
$h_{n}$ & 0.05935281 & 0.06525282 & 0.08533451 & 0.0983018 \\
\hline
\end{tabular*}
\end{table}

%
%
\begin{table}[b]
\tablewidth=150pt
\caption{$h_{0,n}$'s for different $n$'s} \label{ta:parameter}
\begin{tabular*}{\tablewidth}{@{\extracolsep{\fill}}lcc@{}}
\hline
$\bolds n$ & \textbf{10\mbox{,}000} & \textbf{5000} \\
\hline
$h_{0,n}$ & 0.16 & 0.2 \\
\hline
\end{tabular*}
\end{table}

The local bootstrap procedure used here is the same as the one proposed
by Paparoditis and Politis \cite{MR1771485} 
except that here the $Z_i$'s are not lagged variables. For a given
sample $\{ (X_i, Y_i, Z_i) \}_{i=1}^n$, a local bootstrap resample $\{
(X^*_i, Y^*_i, Z^*_i) \}_{i=1}^n$ is generated as follows.
\begin{itemize}
\item Step 1. Draw a random sample $(Z^*_1, \ldots, Z^*_{n})$ from the
empirical cumulative distribution function $\hat{F}_Z$, where
\[
\hat{F}_Z(z) = \frac{1}{n} \sum_{i=1}^n I_{(-\infty, Z_i]}(z).
\]
\item Step 2. For $1 \leq i \leq n$, for each $Z_i^*$ from Step 1, draw
$X_i^*$ and $Y_i^*$ independently from the empirical conditional
cumulative distribution functions $\hat{F}_{X|Z=Z_i^*}$ and $\hat
{F}_{Y|Z=Z_i^*}$, respectively, where
\[
\hat{F}_{X|Z=Z_i^*}(x) = \frac{ \sum_{i=1}^n
k_0((Z_i^*-Z_i)/b)I_{(-\infty, X_i]}(x)}{\sum_{i=1}^n k_0((Z_i^*-Z_i)/b)}
\]
and
\[
\hat{F}_{Y|Z=Z_i^*}(y) = \frac{ \sum_{i=1}^n
k_0((Z_i^*-Z_i)/b)I_{(-\infty, Y_i]}(y)}{\sum_{i=1}^n k_0((Z_i^*-Z_i)/b)}.
\]
The parameters for test 1 with local bootstrap are chosen as follows.
The bandwidth $b$ is taken to be $h_n^{0.4}$, $p_n = q_n = 2$ and
$h_{0,n} = 0.4$, where $h_n$ is as in Table~\ref{ta:hn}.
\end{itemize}

\subsection{Experiments} \label{sec:exp}
The objective of the first experiment is to compare the power of test 1
with that of a Hellinger distance-based test proposed by Su and White
\cite{MR2428851}. 
The critical value for Su and White's test can be determined using the
asymptotic distribution of the test statistic or using local bootstrap.
To distinguish between the two cases, we use test 2A to denote the
asymptotic distribution-based version of Su and White's test and test
2B to denote the local bootstrap version. While test 2B is recommended
by Su and White
\cite{MR2428851}, 
test 2A is used here to save time for computation.

In this experiment, both tests 1 and 2A are carried out for 1000
random samples of size $n=10^4$, where the distribution of $(X,Y,Z)$ is
as in (M2) or (M3). Under~(M2), test 1 is applied to transformed data,
as mentioned in Section \ref{sec:para}. Test 2A is applied to normalized data and
the bandwidth parameter in the kernel estimators in the test statistic
is taken to be $n^{-1/8.5}$, as in
\cite{MR2428851}. 
The power estimates based on data from (M2) and (M3) with $n=10^4$ are
given in Table \ref{ta:M23}. The asymptotic significance level is
$0.05$. It is shown in Table \ref{ta:M23} that power estimates for
test~1 when $a=0$ and $a=0.1$ are larger that those for test 2A.

%
%
\begin{table}[b]
\caption{Power comparison between tests 1 and \textup{2A}} \label{ta:M23}
\begin{tabular*}{\tablewidth}{@{\extracolsep{\fill}}lcccccc@{}}
\hline
& \multicolumn{2}{c}{$\bolds{a=0}$} & \multicolumn{2}{c}{$\bolds{a=0.1}$}
& \multicolumn{2}{c@{}}{$\bolds{a=0.3}$} \\[-4pt]
& \multicolumn{2}{c}{\hrulefill} & \multicolumn{2}{c}{\hrulefill}
& \multicolumn{2}{c@{}}{\hrulefill} \\
& \textbf{Test 1} & \textbf{Test 2A} & \textbf{Test 1}
& \textbf{Test 2A} & \textbf{Test 1} & \textbf{Test 2A} \\
\hline
(M2) & 0.049 & 0.028 & 0.65\phantom{0} & 0.076 & 1 & 0.95 \\
(M3) & 0.041 & 0.029 & 0.572 & 0.119 & 1 & 1\phantom{.00}\\
\hline
\end{tabular*}
\end{table}

To explore the power performance of test 2B without actually running
the local bootstrap procedure, approximate critical values for test 2B
under (M2) and (M3) are used. To obtain these approximate critical
values, note that under (M2) or (M3), for large $n$, a local bootstrap
resample for $a=0.1$ or $a=0.3$ is approximately distributed as a
random sample for the $a=0$ case, so the critical value for test~2B can
be approximated by the 95\% sample quantile of the 1000 test 2A
statistics from the first experiment for the $a=0$ case. Then the power
estimates for test 2B can be approximated by the proportions of the
1000 test 2A statistics from the first experiment under different
alternatives that exceed the approximate critical values. The
approximate power estimates are given in Table \ref{ta:m23_new}. Note
that the approximate power estimates for test 2B are often larger than
the power estimates for test 2A in Table \ref{ta:M23}, which suggests
that test 2B is more powerful than test 2A.

%
%
\begin{table}
\tablewidth=200pt
\caption{Approximated power estimates for test \textup{2B}} \label{ta:m23_new}
\begin{tabular*}{\tablewidth}{@{\extracolsep{\fill}}lcc@{}}
\hline
& $\bolds{a=0.1}$ & $\bolds{a=0.3}$ \\
\hline
(M2) & 0.128 & 0.971 \\
(M3) & 0.241 & 1\phantom{.000}
\\
\hline
\end{tabular*}
\end{table}

To investigate the performance of test 1 when the sample size is
smaller, in the next experiment, power estimates for test 1 are
computed based on 1000 random samples of size $n=5000$ from (M2) and
(M3). The results are given in Table \ref{ta:M23n}. The results for
$n=10^4$ from the first experiment are also included for comparison.
The asymptotic significance level is $0.05$ as before. Table \ref
{ta:M23n} shows that test 1 is more powerful when $n$ is larger.

%
%
\begin{table}[b]
\caption{Test 1 power estimates for $n=5000$ and $n=10^4$} \label{ta:M23n}
\begin{tabular*}{\tablewidth}{@{\extracolsep{\fill}}lcccccc@{}}
\hline
& \multicolumn{2}{c}{$\bolds{a=0}$} & \multicolumn{2}{c}{$\bolds{a=0.1}$}
& \multicolumn{2}{c@{}}{$\bolds{a=0.3}$} \\[-4pt]
& \multicolumn{2}{c}{\hrulefill} & \multicolumn{2}{c}{\hrulefill}
& \multicolumn{2}{c@{}}{\hrulefill} \\
& \textbf{(M2)} & \textbf{(M3)} & \textbf{(M2)} & \textbf{(M3)}
& \textbf{(M2)} & \textbf{(M3)} \\
\hline
$n=5000$ & 0.052 & 0.039 & 0.373 & 0.321 & 0.998 & 1 \\
$n=10^4$ & 0.049 & 0.041 & 0.65\phantom{0} & 0.572 & 1\phantom{.000} & 1
\\
\hline
\end{tabular*}
\end{table}

Finally, for smaller sample size such as $n=500$ or $n=1000$, since the
approximation in Theorem \ref{thm:distribution} does not work well, the
local bootstrap version of test 1 is considered. Here 1000 samples of
size $n$ from (M2) are used, and for each sample, 1000 local bootstrap
resamples are used to determine the rejection region. The level is 0.05.
The power estimates for the test are given in Table \ref{ta:cb}.

In the above results, the power estimates for test 1 are larger when
$a$ is larger. This is expected. Under (M2) or (M3), $E\rho_{p_n,
q_n}^2(Z) = E\rho_{2, 2}^2(Z)$ increases as $a$ increases ($a \in
[0,1]$), so test 1 should be more powerful for larger $a$, if the
approximation in (\ref{eq:consistency2}) and (\ref{eq:uniform}) work.
Table \ref{table7} gives the values of $E\rho_{p_n, q_n}^2(Z)$ for
$a=0.1$ and $0.3$. For (M2), the calculation of $E\rho_{p_n, q_n}^2(Z)$
is done for the transformed $(X,Y,Z)$, which is obtained by applying
the function $\Phi$ to the original $(X,Y,Z)$.

\section{Concluding remarks}
A test statistic for testing conditional independence based on maximal
nonlinear conditional correlation is proposed. Two tests,
tests~1 and 1N, are constructed using the test statistic. Both
tests are consistent and have similar asymptotic properties, as
discussed in Section \ref{sec:test}. Some simulation experiments are
carried out to check the performance of test 1. The simulation results
show that when the sample size $n=10^4$, the power of test 1 is
comparable with that of test 2A.
The simulation results also indicate that test 1 has better power when
$E\rho_{p_n,q_n}^2(Z)$ is larger, as expected.

%
%
\begin{table}
\tablewidth=250pt
\caption{Power estimates for test 1 with local bootstrap} \label{ta:cb}
\begin{tabular*}{\tablewidth}{@{\extracolsep{\fill}}lccc@{}}
\hline
& $\bolds{a = 0}$ & $\bolds{a=0.1}$ & $\bolds{a=0.3}$ \\
\hline
$n=500$ & 0.041 & 0.071 & 0.309 \\
$n=1000$ & 0.033 & 0.099 & 0.531 \\
\hline
\end{tabular*}
\end{table}

%
%
\begin{table}[b]
\tablewidth=250pt
\caption{$E\rho_{p_n, q_n}^2(Z)$ under \textup{(M2)} and \textup{(M3)}}\label{table7}
\begin{tabular*}{\tablewidth}{@{\extracolsep{\fill}}lcc@{}}
\hline
& $\bolds{a=0.1}$ & $\bolds{a=0.3}$ \\
\hline
(M2) & 0.001345575 & 0.01908246 \\
(M3) & 0.002044604 & 0.01765322 \\
\hline
\end{tabular*}
\end{table}

Below are a few remarks.
\begin{enumerate}
\item Equation (\ref{eq:small_nzpq}) requires that $p_n$, $q_n$ and $n_Z$ grow
slowly comparing to $n$. The parameter selection result in Table \ref
{ta:parameter} in Section \ref{sec:sim} seems to agree with such a
requirement. With $n=10^4$, $n_Z$ is only 5 and $p_n=q_n=2$. When
$p_n=q_n=3$, even with $h_{0,n} = 0.4$ (this corresponds to the
smallest $n_Z$ for $n=10^4$), the distribution of the test statistic
cannot be approximated well by the distribution of $\sum_{k=1}^{n_Z}
\lambda_k$.


\item The parameter selection criteria given in Section \ref{sec:sim}
needs to be studied to see whether the asymptotic properties of test 1
still hold using such a criteria.

\item When the distribution of the test statistic cannot be
approximated well by the distribution of $\sum_{k=1}^{n_Z} \lambda_k$,
it is possible to use local bootstrap version of test 1. However, it
takes a lot of time to obtain the bootstrap resamples, so this approach
is recommended when the sample size $n$ is small.

\item In all theorems proved in this paper, it is assumed that the
$(X_i, Y_i, Z_i)$'s are i.i.d. It is also expected that test 1 works
for some stationary weakly dependent data such as the vector ARMA
processes, where the central limit theorem for the i.i.d. case still
applies. However, to carry out the details in the proofs, one needs the
strong approximation result in Lemma \ref{lm:iid_normal}, which is a
stronger result than the central limit theorem and requires a version
of Lemma \ref{lm:iid_ch} that works for dependent data.

\item Test 1 can be modified to work for discrete $Z$. Modification is
necessary since the rate of convergence for each $\hat{\rho}(z_k)$ is
faster in the discrete case.

\item In Lemma \ref{lm:normal} and Theorems
\ref{thm:consistency} and \ref{thm:distribution},
the $z_k$'s are chosen in $\mathcal{Z}(\varepsilon_n)$ so
that they are $\varepsilon_n$-away from the boundary, and it is assumed
that $h_n/\varepsilon_n = O(n^{-\beta})$ to ensure that certain error
terms in the bias/variance calculation are negligible. For
implementation, the condition $h_n/\varepsilon_n = O(n^{-\beta})$ still
leaves some room for choosing $\varepsilon_n$. This problem can be
eliminated by using a kernel function with compact support, as pointed
out by a reviewer. In particular, if the kernel function $k_0$ is
supported on $[-1, 1]^d$, then one can simply take $\varepsilon_n =
h_n$. In such case, even though the condition $h_n/\varepsilon_n =
O(n^{-\beta})$ does not hold, the results in Lemma \ref{lm:normal} and
Theorems \ref{thm:consistency} and \ref{thm:distribution} remain valid.
\end{enumerate}

\section{Proofs}\label{sec:proofs}
\subsection[Proof of Lemma 1]{Proof of Lemma \protect\ref{lm:normal}}
\label{sec:proof_lm_normal}
Recall that for $1 \leq j \leq k_n$,
\[
W_{n, j}(z) = \sqrt{n h_n^d c_K f_Z(z)}  \Biggl(  \Biggl( \sum_{i=1}^{n}
w_i(z) f_{n, j } (X_i, Y_i, z)  \Biggr) - E\bigl(f_{n,j } (X, Y, z)|Z=z\bigr)
\Biggr).
\]
To prove the asymptotic normality of $W_{n,j}(z_k)$'s, we will
approximate $W_{n,j}(z)$ using sums of i.i.d. random variables.
For $1 \leq i \leq n$, let $w_{0,i}(z) = k_0(h_n^{-1}(z-Z_i))$and let
$\hat{f}_Z(z) = n^{-1}h_n^{-d} \sum_{i=1}^{n} w_{0,i}(z)$. Then
$w_i(z) = n^{-1}h_n^{-d} w_{0,i}(z)/\hat{f}_Z(z)$. For $1 \leq j \leq
k_n$, let
\begin{eqnarray*}
\tilde{W}_{n, j}(z) &=& (n h_n^d f_Z(z))^{-1/2} (c_K)^{1/2} \sum
_{i=1}^{n}  \bigl( w_{0,i} (z) f_{n, j } (X_i, Y_i, z)\\
&&\hspace*{121.6pt}{} - Ew_{0,i}(z)
f_{n, j }
(X_i, Y_i, z)  \bigr)
\end{eqnarray*}
and $\tilde{W}_{n, k_n+1} (z)= \sqrt{n h_n^d c_K} (f_Z(z))^{-1/2}
(\hat
{f}_Z(z) - E\hat{f}_Z(z))$, then
\begin{eqnarray*}
W_{n,j} (z)
&= & \frac{f_Z(z)}{\hat{f}_Z(z)}\tilde{W}_{n,j} (z) + \sqrt{nh_n^d c_K
f_Z(z)} E\bigl(f_{n,j}(X,Y,z)|Z=z\bigr)  \biggl(
\frac{f_Z(z)}{\hat{f}_Z(z)} - 1  \biggr) \\
&&{} + \frac{\sqrt{nh_n^d
c_Kf_Z(z)}}{\hat{f}_Z(z)}  \bigl(
h_n^{-d} E(w_{0,1}(z)f_{n,j}(X_1,Y_1,z))\\
&&\hspace*{80.4pt}{} - E\bigl(f_{n,j}(X,Y,z)|Z=z\bigr) f_Z(z)
\bigr) \\
& = & \hat{W}_{n,j}(z) + \sum_{\ell=1}^4 R_{\ell,n,j}(z),
\end{eqnarray*}
where
$\hat{W}_{n,j}(z) = \tilde{W}_{n,j} (z) - \tilde{W}_{n,k_n+1}
(z)E(f_{n,j}(X,Y,z)|Z=z)$,
\begin{eqnarray*}
R_{1,n,j} (z) &=& \biggl( \frac{f_Z(z)}{\hat{f}_Z(z)} -1  \biggr) \tilde
{W}_{n,j} (z),
\\
R_{2,n,j} (z) &=& \frac{\sqrt{n
h_n^d c_Kf_Z(z)}}{\hat{f}_Z(z)}  \bigl(
h_n^{-d}E(w_{0,1}(z)f_{n,j}(X_1,Y_1,z))\\
&&\hspace*{67.5pt}{} - E\bigl(f_{n,j}(X,Y,z)|Z=z\bigr) f_Z(z) \bigr),
\\
R_{3,n,j} (z)
&=& \frac{ \sqrt{nh_n^d c_K} E(f_{n,j}(X,Y,z)|Z=z) (f_Z(z) - \hat
{f}_Z(z))^2 }{\hat{f}_Z(z) \sqrt{f_Z(z)} }
\end{eqnarray*}
and
\[
R_{4,n,j} (z)= -
\frac{\sqrt{nh_n^d c_K}}{\sqrt{f_Z(z)}}E\bigl(f_{n,j}(X,Y,z)|Z=z\bigr)
\bigl(E\hat{f}_Z(z) - f_Z(z)  \bigr).
\]
We will complete the proof by showing that the following results hold
for $T_n=\exp( -(\ln n)^{1/9})$.
\begin{enumerate}[(C3)]
\item[(C1)] $\sum_{j=1}^{k_n} \sum_{k=1}^{n_Z}  ( \sum_{\ell=1}^4
R_{\ell,n,j}(z_k)  )^2=O_p(T_n)$.

\item[(C2)] There exist random variables $N_{1, j, k}$ and
$\varepsilon
_{1,j,k}\dvtx 1\leq j \leq k_n$, $1\leq k \leq n_Z$ such that
the joint distribution of $(N_{1, j, k}+\varepsilon_{1,j,k})_{j,k}$ is
the same as that of $( \hat{W}_{n,j} (z_k)
)_{j,k}$, $N_{1, j, k}$'s are jointly normal with $EN_{1, j, k}=0$ and
$\operatorname{Cov}(N_{1, j, k}, N_{1, \ell, k^*}) = \operatorname{Cov}(\hat{W}_{n,j}
(z_k), \hat{W}_{n,\ell} (z_{k^*}))$ and $\sum_{j=1}^{k_n} \sum
_{k=1}^{n_Z} \varepsilon_{1,j,k}^2 = O_p(T_n)$.

\item[(C3)] There exist random variables $N_{2, j, k}$ and
$\varepsilon
_{2,j,k}\dvtx 1\leq j \leq k_n$, $1\leq k \leq n_Z$ such
that the joint distribution of $(N_{2, j, k}+\varepsilon
_{2,j,k})_{j,k}$ is the same as that of
$(N_{1, j, k})_{j,k}$, $N_{2, j, k}$'s are jointly normal with $EN_{2,
j, k}=0$ and
\begin{eqnarray*}
&&
\operatorname{Cov}(N_{2, j, k}, N_{2, \ell, k^*}) \\
&&\qquad =
\cases{\operatorname{Cov}\bigl(f_{n,j } (X, Y,z_k), f_{n, \ell} (X,
Y,z_k)|Z=z_k\bigr), &\quad
if $k = k^*$; \cr
0, &\quad otherwise,}
\end{eqnarray*}
and $\sum_{j=1}^{k_n} \sum_{k=1}^{n_Z}
\varepsilon_{2,j,k}^2 = O_p(T_n)$.
\end{enumerate}
Note that Lemma \ref{lm:normal} follows from (C1)--(C3) since one can
construct random variables $\tilde{N}_{2,j,k}$,
$\tilde{\varepsilon}_{2,j,k}$, $\tilde{\varepsilon}_{1,j,k}$ and
$R_{5, n,j,k}\dvtx 1\leq j \leq k_n$, $1\leq k \leq n_Z$ on the
same probability space such that the joint distribution of $(\tilde
{N}_{2,j,k}, \tilde{\varepsilon}_{2,j,k})_{j,k}$ is the same
as that of $(N_{2, j, k}, \varepsilon_{2,j,k})_{j,k}$, the joint
distribution of $(\tilde{\varepsilon}_{1,j,k},
\tilde{N}_{2,j,k}+\tilde{\varepsilon}_{2,j,k})_{j,k}$ is the same as
that of $(\varepsilon_{1,j,k}, N_{1, j, k}
)_{j,k}$, and the joint distribution of $(R_{5, n,j,k}, \tilde
{N}_{2,j,k}+\tilde{\varepsilon}_{2,j,k}+\tilde{\varepsilon}_{1,j,k})_{j,k
}$ is the same as that of $(\sum_{\ell=1}^4 R_{\ell,n,j}(z_k),\break \hat
{W}_{n,j}(z_k))_{j,k}$. Take $W_{n,1,j,k}=\tilde{N}_{2,j,k}$
and $W_{n,2,j,k} = \tilde{\varepsilon}_{2,j,k}+\tilde{\varepsilon
}_{1,j,k}+R_{5, n,j,k}$, then we have Lemma~\ref{lm:normal}.

To establish (C1)--(C3), we need certain expectations and covariances,
which are computed below.
Under (R1)--(R3) and the conditions that\break $\int u k_0(u) \,du = 0$ and
$\sigma_0^2 = \int\| u \|^2 k_0(u) \,du < \infty$, for $z \in\mathcal{Z}
(\varepsilon_n)$, we have
%
%
\begin{eqnarray} \label{eq:bias}
&& (h_n^d)^{-1} E ( w_{0,1} (z) f_{n, j } (X_1,
Y_1, z)  )
\nonumber\\[-8pt]\\[-8pt]
&&\qquad = E \bigl(f_{n, j } (X, Y,z)|Z=z\bigr)f_Z(z) + r_{n,j,1}(z) C_n
h_n^2,\nonumber
\end{eqnarray}
where
\begin{eqnarray*}
r_{n,j,1}(z) &=& c_0 \int h(x,y) \,d\mu(x,y)\\
&&{}  \times\bigl( 2 d \sigma_0^2
\theta
_{n,j,1} + \theta_{n,j,2} h_n^{-2} (2+h_n) \gamma_4^d \exp (
-\gamma_5 \varepsilon_n^2 h_n^{-2}  )  \bigr),
\end{eqnarray*}
$|\theta_{n,j,1}|$, $|\theta_{n,j,2}| \leq1$, and $\gamma_4$ and
$\gamma_5$ are positive constants that depend on $\gamma_2$
and $\gamma_3$ only. Also, for $k \neq k^*$, $z_k$, $z_k^* \in
\mathcal{Z}
(\varepsilon_n)$, we have
%
%
\begin{eqnarray}\label{eq:cov1}\quad
&& (h_n^d)^{-2} \operatorname{Cov}  ( w_{0,1}(z_k) f_{n, j } (X_1, Y_1, z_k),
w_{0,1}(z_{k^*}) f_{n, \ell} (X_1, Y_1, z_{k^*})  )
\nonumber\\
&&\qquad = \theta_{j,\ell, k, k^*} (h_n^d)^{-2} (\gamma
_2)^{2d} \exp( -0.5\gamma_3 h_n^{-2} \| z_k -
z_{k^*} \|^2) C_n^2 \nonumber\\
&&\qquad\quad{} - f_Z(z_k) f_Z(z_{k^*}) E \bigl(f_{n, j
} (X, Y,z_k)|Z=z_k\bigr) E \bigl(f_{n, \ell} (X, Y,z_{k^*})|
Z=z_{k^*}\bigr) \nonumber\\[-8pt]\\[-8pt]
&&\qquad\quad{} - f_Z(z_k)E \bigl(f_{n, j } (X, Y,z_k)|Z=z_k\bigr)
r_{n,\ell,1}(z_{k^*}) C_n h_n^2 \nonumber\\
&&\qquad\quad{} -f_Z(z_{k^*})E \bigl(f_{n, \ell} (X,
Y,z_{k^*})|Z=z_{k^*}\bigr)r_{n,j,1}(z_k) C_n h_n^2 \nonumber\\
&&\qquad\quad{} - r_{n,j,1}(z_k)
r_{n,\ell,1}(z_{k^*}) C_n^2 h_n^4,\nonumber
\end{eqnarray}
where $|\theta_{j,\ell, k, k^*}| \leq1$.
Finally, for $z\in\mathcal{Z}(\varepsilon_n)$,
%
%
\begin{eqnarray} \label{eq:cov2}
&& (h_n^d)^{-1} \operatorname{Cov}  ( w_{0,1}(z) f_{n, j } (X_1, Y_1, z),
w_{0,1}(z) f_{n, \ell} (X_1, Y_1, z)  ) \nonumber\\
&&\qquad = f_Z(z) E\bigl(f_{n, j } (X, Y,z) f_{n,\ell}(X,Y,z)|Z=z\bigr) \int
k_0^2(u)\,
du + r_{n,j, \ell, 2}(z) C_n^2 h_n \nonumber\\
&&\qquad\quad{} - h_n^d f_Z^2(z)E \bigl(f_{n, j } (X, Y,z)|Z=z\bigr) E \bigl(f_{n, \ell} (X,
Y,z)|Z=z\bigr) \nonumber\\[-8pt]\\[-8pt]
&&\qquad\quad{} - h_n^{d+2} C_n r_{n,j,1}(z)
f_Z(z) E \bigl(f_{n, \ell} (X, Y,z)|Z=z\bigr) \nonumber\\
&&\qquad\quad{} - h_n^{d+2} C_n
r_{n, \ell,1}(z) f_Z(z) E \bigl(f_{n, j } (X, Y,z)|Z=z\bigr)
\nonumber\\
&&\qquad\quad{} -h_n^{d+4} C_n^2 r_{n,j,1}(z) r_{n, \ell,1}(z)\nonumber
\end{eqnarray}
and
%
%
\begin{equation} \label{eq:moment3} h_n^{-d} E (
w_{0,1}(z)f_{n,j}(X_1,Y_1,z)  )^3 \leq C_n^3 c_0 \int k_0^3(u) \,du,
\end{equation}
where
\[
| r_{n,j,\ell, 2}(z) | \leq2 c_0 \int h(x,y)\, d\mu(x,y)  \biggl( \sqrt
{d} \int\| u \| k_0^2(u) \,du + h_n^{-1} \gamma_6^d
e^{-\gamma_7 \varepsilon_n^2/h_n^2}  \biggr)
\]
for some positive constants $\gamma_6$ and $\gamma_7$ that
depend on $\gamma_2$ and $\gamma_3$ only. Below we will prove (C1)--(C3).
\begin{pf*}{Proof of (C1)}
Let $S_n=\sum_{k=1}^{n_z} (\hat{f}_Z(z_k) -
f_Z(z_k))^2$ and $A_n = \{ \sqrt{S_n} < \min\{1$, $(2c_1)^{-1} \} \}$.
From (\ref{eq:bias}) and (\ref{eq:cov2}), $ES_n=O(n_Z(h_n^4+
(nh_n^d)^{-1}) ) = O(n_Z(n\times h_n^d)^{-1})$ and $1/f_Z(z_k) \leq c_1$ for
all $k$,
$P(A_n^c) \rightarrow0$ as $n\rightarrow\infty$. From (\ref
{eq:bias}), on $A_n$,
\begin{eqnarray*}
&& \sum_{j=1}^{k_n} \sum_{k=1}^{n_Z}  \Biggl( \sum_{\ell=1}^4
|R_{\ell
,n,j}(z_k)|  \Biggr)^2 \\
&&\qquad \leq O(1)  \Biggl( S_n  \Biggl( \sum_{j=1}^{k_n}
\sum_{k=1}^{n_Z} \tilde{W}_{n, j}^2(z_k)  \Biggr) + k_n n_Z C_n^2 (n
h_n^{d+4}) + k_n C_n^2 n h_n^d S_n^2  \Biggr), \\
\end{eqnarray*}
and it follows from (\ref{eq:cov2}) that
\[
E  \Biggl( \sum_{j=1}^{k_n} \sum_{k=1}^{n_Z} \tilde{W}_{n, j}^2
(z_k)  \Biggr) = O(k_n n_Z C_n^2).
\]
Take
\[
T_{1,n} = \frac{k_n n_Z^2 C_n^2}{nh_n^d} + k_n n_Z C_n^2 n h_n^{d+4},
\]
then (C1) holds with
$T_n =\exp( -(\ln n)^{1/9}) $ since $T_{1,n}= O(T_n)$.
\end{pf*}

The proof of (C2) is based on the following lemma, which deals
with the normal approximation of sum of i.i.d.
random vectors.
\begin{lm} \label{lm:iid_normal}
Suppose that $X_1, \ldots, X_n$ are i.i.d. random vectors in
$R^{d_1}$ with mean 0 and variance $\Sigma$.
Suppose that there exist positive constants $C$, $a_2$ and $a_3$ such
that $1 \leq a_2 \leq a_3 \leq C$,
$\| X_1 \| \leq C$ and $E\| X_1 \|^k \leq a_k^k$ for $k=2$, 3. Then for
$T\geq1$, there exist random vectors $S$ and
$Y$ on the same probability space such that $S$ is distributed as $(X_1
+ \cdots+ X_n)/\sqrt{n}$, $Y$ is multivariate
normal with mean 0 and variance $\Sigma$ and for $n \geq(25/(16a_2^2)
+ 25d_1/12)C^2T^4 \exp(3T^2/16)$,
\[
P(\| S-Y \| \geq\alpha) \leq\alpha,
\]
if
\[
\alpha\geq\frac{33.75a_3^3}{\sqrt{n}} (12)^{d_1} e^{(d_1+3)T^2/8} +
(48)^{d_1} e^{-3T^2/(32a_2^2)}.
\]
\end{lm}

The proof of Lemma \ref{lm:iid_normal} is given in Section
\ref{sec:iid_normal}. To prove (C2), note that
$\tilde{W}_{n,j}(z_k) = \sum_{i=1}^n (g_{n,j,k}(X_i, Y_i, Z_i) -
Eg_{n,j,k}(X_i, Y_i, Z_i))/\sqrt{n}$, where
\begin{eqnarray*}
&& g_{n,j,k}(X_i, Y_i, Z_i) \\
&&\qquad = \frac{\sqrt{c_K}}{\sqrt
{f_Z(z_k)h_n^d}} k_0 \biggl( \frac{z_k - Z_i}{h_n}  \biggr)\\
&&\qquad\quad{}\times
\bigl( f_{n,j}(X_i, Y_i, z_k) - E\bigl(f_{n,j}(X, Y, z_k) |Z=z_k\bigr)  \bigr).
\end{eqnarray*}
From (\ref{eq:bias})--(\ref{eq:moment3}),
we have
\begin{eqnarray*}
\Biggl( \sum_{j=1}^{k_n} \sum_{k=1}^{n_Z}  \bigl( g_{n,j,k}(X_i, Y_i,
Z_i) - Eg_{n,j,k}(X_i, Y_i, Z_i)
 \bigr)^2  \Biggr)^{1/2} & \leq & \frac{O(1)C_n\sqrt{k_n n_Z}}{\sqrt
{h_n^d}},
\\
\Biggl( \sum_{j=1}^{k_n} \sum_{k=1}^{n_Z}
E \bigl( g_{n,j,k}(X_i, Y_i, Z_i) - Eg_{n,j,k}(X_i, Y_i, Z_i)
\bigr)^2 \Biggr)^{1/2} &\leq& O(1) C_n \sqrt{k_n n_Z}
\end{eqnarray*}
and
\begin{eqnarray*}
&&  \Biggl( E  \Biggl( \sum_{j=1}^{k_n} \sum_{k=1}^{n_Z}
\bigl( g_{n,j,k}(X_i, Y_i, Z_i) -
Eg_{n,j,k}(X_i, Y_i, Z_i)  \bigr)^2  \Biggr)^{3/2}  \Biggr)^{1/3} \\
&&\qquad
\leq C_n \sqrt{k_n n_Z} h_n^{-d/6} O(1).
\end{eqnarray*}
Note that
for every constant $M>0$, the condition
\[
n \geq \biggl( \frac{25}{16} + \frac{25k_nn_Z}{12}  \biggr)
\biggl(\frac{MC_n
\sqrt{k_n n_Z}}{\sqrt{h_n^d}}  \biggr)^2 T_{3,n}^4 e^{3T_{3,n}^2/16}
\]
holds for large $n$ with $T_{3,n} = (\ln n)^{1/8}$,
so Lemma \ref{lm:iid_normal} is applicable. From Lem\-ma~\ref
{lm:iid_normal}, (C2) holds with any $T_n$ such that $T_{2,n} =
O(T_n)$, where
\[
T_{2,n} = \frac{(C_n \sqrt{k_n n_Z})^612^{2k_nn_Z}
e^{(k_nn_Z+3)T_{3,n}^2/4}}{nh_n^d} + (48)^{2k_nn_Z}
e^{-\gamma T_{3,n}^2/(C_n \sqrt{k_nn_Z})^2},
\]
$\gamma>0$ is a constant. Since
$T_{2,n} = O( \exp( -\gamma_1 (\ln n)^{1/8} ))$ for some constant
\mbox{$\gamma_1 > 0$}, (C2) holds with $T_n=\exp( -(\ln n)^{1/9})$.

The proof of (C3) is based on the following result.
\begin{fa} \label{fa:sqrt} Suppose that $A$ and $B$ are $d_1 \times
d_1$ nonnegative definite matrices. Then
\[
\bigl\| \sqrt{A} - \sqrt{B} \bigr\| \leq d_1^{3/4} \sqrt{ \| A - B\| }.
\]
\end{fa}

The proof of Fact \ref{fa:sqrt} is given at the end of the proof of
(C3). Note that Fact~\ref{fa:sqrt} implies the following: suppose that
$X_0$ and $Y_0$ are two $d_1 \times1$ normal vectors of mean 0 and
covariance matrices $A$ and $B$, respectively.
Let $Z$ be a $d_1 \times1$ normal vector whose elements are i.i.d.
$N(0,1)$. Then $\sqrt{A} Z$ is distributed as $X_0$
and $\sqrt{B} Z$ is distributed as $Y_0$ and
\begin{eqnarray*}
\bigl\| \sqrt{A} Z - \sqrt{B} Z \bigr\|^2 &\leq& \bigl\| \sqrt{A} - \sqrt{B} \bigr\|^2
\| Z \|
^2 \leq d_1^{3/2}
\| A - B \| \| Z \|^2 \\
&=& O_p(d_1^{5/2} \| A-B \|).
\end{eqnarray*}
Therefore, (C3) holds if $\operatorname{Cov}(\hat{W}_{n,j}(z_k), \hat{W}_{n,\ell}(z_{k^*}))$
is close to
\[
\operatorname{Cov}\bigl(f_{n,j } (X, Y,z_k), f_{n, \ell} (X, Y,z_k)|Z=z_k\bigr) \delta_{k,k^*},
\]
where $ \delta_{k, k^*}$ is 1 if $k=k^*$ and
is 0 otherwise. From (\ref{eq:bias})--(\ref{eq:moment3}), we have
\begin{eqnarray*}
&& \sum_{j,\ell, k, k^*}  \bigl( \operatorname{Cov}(\hat{W}_{n,j}(z_k), \hat
{W}_{n,\ell
}(z_{k^*})) \\
&&\qquad\quad\hspace*{-3.1pt}{} -
\operatorname{Cov}\bigl(f_{n,j } (X, Y,z_k), f_{n, \ell} (X, Y,z_k)|Z=z_k\bigr)
\delta_{k,k^*}
\bigr)^2 \\
&&\qquad =h_n C_n^2 (k_n n_Z)^2 O(1),
\end{eqnarray*}
so (C3) holds with $T_n =\exp( -(\ln n)^{1/9})$ since $(k_n
n_Z)^{5/2}\sqrt{h_n C_n^2 (k_n n_Z)^2} = O(\exp( -(\ln n)^{1/9}))$.
%
%
\begin{pf*}{Proof of Fact \ref{fa:sqrt}}
Consider first the case where $A$ is diagonal. Let $D$ be a diagonal
matrix such that
$B=Q^TDQ$ for some $Q$ such that $QQ^T=I$.
Let $D = \operatorname{diag} (\lambda_1, \ldots, \lambda_{d_1})$, $A =
\operatorname{diag} (\alpha_1, \ldots, \alpha_{d_1})$, $Q = (q_{i,j})$
and $E = B-A= (e_{i,j})$. Let $q_i$ be the $i$th column of $Q$, then
$q_i^TDq_j = \alpha_i \delta_{i,j} + e_{i,j}$, where
$\delta_{i,j} = 1$ for $i=j$ and $\delta_{i,j} = 0$, otherwise. Write
$Dq_k = \sum_{j=1}^{d_1} (q_k^T D q_j) q_j$, then
%
\begin{eqnarray*}
\bigl\| \sqrt{D}q_k - \sqrt{\alpha_k} q_k \bigr\|^2 &=& \sum_{j=1}^{d_1}
\bigl(\sqrt{\lambda_j} q_{j,k} -\sqrt{\alpha_k} q_{j,k}\bigr)^2 \\
&=& \sum_{j=1}^{d_1}  \bigl( \sqrt{\lambda_j |q_{j,k}|} -\sqrt{
\alpha
_k|q_{j,k}| }  \bigr)^2 |q_{j,k}| \\
&\leq&\sum_{j=1}^{d_1}
\bigl| \lambda_j |q_{j,k}| - \alpha_k|q_{j,k}|  \bigr| |q_{j,k}| \\
&\leq& \Biggl( \sum_{j=1}^{d_1} ( \lambda_j q_{j,k} - \alpha_k
q_{j,k} )^2  \Biggr)^{1/2}  \Biggl( \sum_{j=1}^{d_1} q_{j,k}^2
\Biggr)^{1/2} \\
&=&  \Biggl( \sum_{j=1}^{d_1} e_{k,j}^2  \Biggr)^{1/2}
\end{eqnarray*}
and
\begin{eqnarray*}
\bigl\| \sqrt{Q^TD Q} - \sqrt{A} \bigr\|^2 &=& \sum_{i=1}^{d_1} \sum
_{j=1}^{d_1}  \bigl( q_i^T \sqrt{D} q_j - q_i^T
\sqrt{\alpha_j} q_j  \bigr)^2 \\
&\leq& \sum_{i=1}^{d_1} \sum
_{j=1}^{d_1} \bigl\| \sqrt{D} q_j -\sqrt{\alpha_j} q_j \bigr\|^2
\\
&\leq& d_1 \sum_{j=1}^{d_1}  \Biggl( \sum_{\ell=1}^{d_1} e_{j, \ell}^2
 \Biggr)^{1/2} \\
& \leq& (d_1)^{3/2}  \Biggl( \sum_{j=1}^{d_1}
\sum_{\ell=1}^{d_1} e_{j, \ell}^2  \Biggr)^{1/2},
\end{eqnarray*}
so the result in Fact \ref{fa:sqrt} holds if $A$ (or $B$) is diagonal.
For general $A$ and $B$, write $A=P^T A_0 P$ and $B =Q^TDQ$, where
$A_0$ and $D$ are diagonal and
$P^TP=Q^TQ=I$. Let $B_0= PQ^TDQP^T$, then we have
\begin{eqnarray*}
\bigl\| \sqrt{A} - \sqrt{B} \bigr\| & = & \bigl\| P^T\sqrt{A_0} P - Q^T\sqrt{D} Q
\bigr\| \\
& = & \bigl\| \sqrt{A_0} - PQ^T\sqrt{D} QP^T \bigr\| \leq d_1^{3/4} \sqrt{ \|
A_0 - B_0 \| } \\
& = & d_1^{3/4} \sqrt{ \| P^TA_0P - P^TB_0P \| } = d_1^{3/4} \sqrt{ \|
A - B \| }.
\end{eqnarray*}
The proofs of Fact \ref{fa:sqrt} and Lemma \ref{lm:normal} are complete.
\end{pf*}

\subsubsection[Proof of Lemma 2]{Proof of Lemma \protect\ref{lm:iid_normal}}
\label{sec:iid_normal}

The proof Lemma \ref{lm:iid_normal} is based on several facts, which
are taken directly or adapted from some existing results
and are stated/proved below in Lemmas \ref{lm:lm2.1_bp}--\ref{lm:iid_ch}.

In the statements of Lemmas \ref{lm:lm2.1_bp} and \ref{lm:ch_prohorov},
$(S_0, d_0)$ is a metric space, $\mathcal{B}$ denotes the collection
of Borel
sets in $(S_0, d_0)$, and for two measures $\mu_1$ and $\mu_2$
defined on $\mathcal{B}$, $\rho_0( \mu_1, \mu_2)$ denotes the
Prohorov distance
of $\mu_1$ and $\mu_2$, which is defined as
\[
\rho_0( \mu_1, \mu_2) = \inf\{ \epsilon> 0\dvtx \mu_1(A) < \mu_2
(A^{\epsilon} )+ \epsilon\mbox{,  for all } A \in\mathcal{B}\},
\]
where $A^{\epsilon} =
\{ x\dvtx d^*(x, A) < \epsilon\}$ and $d^*(x, A) = \inf\{ d_0(x, y)\dvtx y
\in A \}$. Here are Lemmas \ref{lm:lm2.1_bp}--\ref{lm:iid_ch}.
\begin{lm}[(Lemma 2.1 in Berkes and Philipp \cite{MR515811})]
\label{lm:lm2.1_bp}  
Suppose that $P_1$ and $P_2$ are
two measures defined on $\mathcal{B}$ and $\rho_0 ( P_1, P_2) <
\alpha$. Then
there exists a probability measure
$Q$ on the Borel sets of $S_0 \times S_0$ with marginals $P_1$ and
$P_2$ such that
\[
Q\{ (x,y)\dvtx d_0(x,y) > \alpha\}
\leq\alpha.
\]
\end{lm}
\begin{lm}[(Adapted from Lemma 2.2 in \cite{MR515811})] \label{lm:ch_prohorov}
Suppose that\vspace*{1pt} $F$ and $G$ are two distributions on $R^{d_1}$ with
characteristic functions $f$ and
$g$, respectively. Then for $\sigma\in(0, 1]$ and $T > 0$, the
Prohorov distance $\rho_0(F,G) \leq\alpha$, where
\begin{eqnarray*}
\alpha &=& \sigma T + 3 (2^{d_1}) e^{-{3T^2}/{32}} + \biggl( \frac
{T}{\pi}  \biggr)^{d_1} \int|f(u) - g(u)| e^{- { \sigma^2 \| u
\|
^2}/{2}} \,du\\
&&{}
+ F \biggl(  \biggl\{ x \dvtx \| x \|
\geq\frac{T}{2}  \biggr\}  \biggr).
\end{eqnarray*}
\end{lm}
\begin{pf}
Let $H$ be the $N(0, \sigma^2 I )$
distribution on $R^{d_1}$, where $I$ is the identity
matrix and $\sigma> 0$. Let $F_1$ be the convolution of $F$ and $H$
and $G_1$ be the convolution of $G$ and $H$.
Then
%
%
\begin{equation} \label{eq:prohorov_1} \rho_0(F,G) \leq\rho_0(F_1,
G_1) + 2 \max\bigl\{ r, H(\{ x\dvtx \| x \| \geq r \}) \bigr\}\qquad
\mbox{for every $r>0$.}\hspace*{-32pt}
\end{equation}
Let $f_1$, $g_1$ and $h$ be the characteristic functions of $F_1$,
$G_1$ and $H$, respectively, and let
$\gamma_F$ and $\gamma_G$ be the densities of $F_1$ and $G_1$,
respectively. Then
\begin{eqnarray*}
|\gamma_F(x) - \gamma_G(x)| & = &
(2\pi)^{-d_{1}}  \biggl| \int e^{-i u^T x} \bigl(f_1(u)-g_1(u)\bigr) \,du  \biggr|
\\
& \leq& (2\pi)^{-d_{1}} \int|f(u) - g(u)| |h(u)| \,du,
\end{eqnarray*}
which implies that for every borel set $B$ in $R^{d_1}$,
\begin{eqnarray*}
&& F_1(B) - G_1(B) \\
&&\qquad\leq F_1(B \cap\{ x\dvtx \|x\| \leq T \}) -
G_1(B \cap\{ x\dvtx \|x\| \leq T \}) + F_1(\{ x\dvtx \|x\| \geq T \}) \\
&&\qquad\leq {\int_{ \{ x\dvtx \| x \| \leq T \} }} | \gamma_F(x) -
\gamma_G(x)|\,dx + F(\{ x\dvtx \|x\| \geq T/2 \})\\
&&\qquad\quad{} + H(\{ x\dvtx \|x\| \geq T/2 \}) \\
&&\qquad\leq \fontsize{10.2}{12}\selectfont{\underbrace{  \biggl( \frac{T}{\pi}  \biggr)^{d_1} \int
|f(u) - g(u)| |h(u)| \,du + F(\{ x\dvtx \|x\| \geq T/2 \}) + H(\{ x\dvtx \|x\| \geq T/2 \}
) }_{\mathit{II}}}.
\end{eqnarray*}
Note that $\mathit{II}$ is an upper bound for the Prohorov distance $\rho_0(F_1,
G_1)$, so for $r \leq T/2$, it follows from (\ref{eq:prohorov_1}) that
\begin{eqnarray*}
\rho_0(F,G) & \leq& \mathit{II} + 2 r + 2H(\{ x\dvtx \| x \| \geq r \}) \\
& \leq&  \biggl( \frac{T}{\pi}  \biggr)^{d_1} \int|f(u) - g(u)| |h(u)|
\,du + F(\{ x\dvtx \|x\| \geq T/2 \}) + 2r \\
&&{} + 3 P\bigl( \chi^2(d_1) \geq(r/\sigma)^2 \bigr).
\end{eqnarray*}
Since $h(u) = e^{-\sigma^2 \| u \|^2/2}$ and
%
%
\begin{eqnarray} \label{eq:chisq}
P\bigl( \chi^2(d_1) \geq A \bigr)
&\leq&
e^{-t A} Ee^{t \chi^2(d_1)}
 |_{t=3/8} \nonumber\\[-8pt]\\[-8pt]
&=& e^{-3A/8} (2^{d_1}) \qquad\mbox{for every $A > 0$}.\nonumber
\end{eqnarray}
Lemma \ref{lm:ch_prohorov} holds if $r=\sigma T/2$ and
$\sigma\in(0, 1]$.
\end{pf}
\begin{lm}[(Adapted from Theorem 1(a) in pages
204--208 in Gnedenko and Kolmogorov \cite{gnedenkolimit1968})] \label{lm:iid_ch}
Suppose that $X_1, \ldots, X_n$ are i.i.d. random vectors with mean
0 and variance $\Sigma$. Suppose that $C$ and $a$
are positive constants such that $\| X_1 \| \leq C$, $a \leq C$ and
$E\| X_1\|^k \leq a^k$ for $k=2$, 3. Let $f_n$ be the characteristic
function of $(X_1 + \cdots+ X_n)/\sqrt{n}$. Then
\[
\biggl| f_n(u) - \exp \biggl( -\frac{1}{2} u^T \Sigma u  \biggr)
\biggr|
\leq\frac{ 0.25\| u \|^3 a^3 }{
\sqrt{n} },
\]
if $\| u \| \leq(0.4\sqrt{n})/C$.
\end{lm}
\begin{pf}
Consider first the case where $X_1$ is
univariate. Let $U = f_1(u/\sqrt{n}) -1$, then
\[
U = \frac{\theta^*_1 EX_1^2}{2}  \biggl( \frac{u}{\sqrt{n}}  \biggr)^2
\]
and
\[
U = \frac{EX_1^2}{2}  \biggl( \frac{iu}{\sqrt{n}}  \biggr)^2 + \frac{
\theta_1E|X_1|^3}{3!}  \biggl( \frac{u}{\sqrt{n}}  \biggr)^3,
\]
where $|\theta^*_1| \leq1$ and $|\theta_1| \leq1$.
Suppose that $|u| \leq(0.4\sqrt{n})/C$, then $|U| < 0.1$ and
\[
\log(1+ U) =U + 0.62 \theta_2 U^2 ,
\]
where
$|\theta_2| \leq1$. Let $V = \log f_n(u) +E(X_1^2)u^2/2 =
E(X_1^2)u^2/2+ n \log(1+ U)$, then
\begin{eqnarray*}
V & = & \frac{ n\theta_1E|X_1|^3u^3}{3!n^{3/2}} + (0.62) n \theta_2
 \biggl( \frac{EX_1^2}{2}  \biggl( \frac{iu}{\sqrt{n}}  \biggr)^2
+ \frac{ \theta_1E|X_1|^3}{3!}  \biggl( \frac{u}{\sqrt{n}}  \biggr)^3
 \biggr)^2 \\
& = & \frac{\lambda_1 |u|^3 a^3 }{6\sqrt{n} } + 0.62  \biggl( \frac
{\lambda_2 a^4 u^4}{4n} + \frac{\lambda_3 a^5 |u|^5}{6 (\sqrt{n})^3}
+ \frac{\lambda_4 a^6 u^6 }{36n^2}  \biggr) \\
& = & \frac{|u|^3 a^3 }{\sqrt{n} }  \biggl( \frac{\lambda_1}{6} + 0.62
 \biggl( \frac{\lambda_2 a|u|}{4\sqrt{n}} + \frac{\lambda_3
a^2u^2}{6 n }
+ \frac{\lambda_4 a^3 |u|^3}{36 (\sqrt{n})^3 }  \biggr)  \biggr),
\end{eqnarray*}
where $|\lambda_k| \leq1$ for $k=1$, 2, 3, 4. Since $a|u|/\sqrt{n}
\leq0.4$,
\[
V = \frac{\theta_3 (0.25) |u|^3 a^3 }{\sqrt{n} },
\]
where $|\theta_3| \leq1$. Since $e^V = 1 + \theta_4 |V| e^{|V|}$,
where $|\theta_4| \leq1$,
\begin{eqnarray*}
f_n(u) & = & \exp \biggl( - \frac{E(X_1^2)u^2}{2}\biggr)
\bigl(1 + \theta_4 |V| e^{|V|} \bigr) \\
& = & \exp \biggl( - \frac
{E(X_1^2)u^2}{2}  \biggr) + \theta_5  \biggl( \frac{ 0.25|u|^3a^3 }{
\sqrt{n} }
 \biggr) e^{|V| - E(X_1^2)u^2/2},
\end{eqnarray*}
where $|\theta_5| \leq1$. To find an upper bound for $|V| -
E(X_1^2)u^2/2$, note that
\[
\biggl| n U + \frac{E(X_1^2)u^2}{2}\biggr| = \frac{|\theta_1| E|X_1|^3
|u|^3}{6 \sqrt{n}} \leq
\frac{ C EX_1^2 |u|^3}{6 \sqrt{n}} \leq\frac{(0.4) u^2 E(X_1^2)}{6},
\]
$n|U| = |\theta^*_1| u^2 E(X_1^2)/2 \leq u^2 E(X_1^2)/2$
and
\[
\bigl|n \bigl(\log(1+U) - U\bigr)\bigr| = 0.62 n |\theta_2 U^2| \leq0.62 (0.1)  \biggl(
\frac{ E(X_1^2)u^2}{2}  \biggr)
\]
since $|U| < 0.1$. Therefore,
\begin{eqnarray*}
|V| -\frac{u^2E(X_1^2)}{2} & = &  \biggl| \frac{E(X_1^2)u^2}{2} + n U +
n \bigl(\log(1+U) - U\bigr)  \biggr| -\frac{u^2E(X_1^2)}{2} \\
& \leq& \frac{(0.4) u^2 E(X_1^2)}{6} + \frac{0.062E(X_1^2)u^2 }{2}
-\frac{u^2E(X_1^2)}{2} \leq0
\end{eqnarray*}
and
Lemma \ref{lm:iid_ch} holds for the univariate case. The result for
the general case
can be obtained by applying the univariate result with $u$ and $X_i$
replaced by $\| u \|$ and $Y_i = u^T X_i/\| u \|$.
\end{pf}

Now we are ready to prove Lemma \ref{lm:iid_normal}.
\begin{pf*}{Proof of Lemma \ref{lm:iid_normal}}
Let $f_n$ be the
characteristic function of $(X_1 + \cdots+ X_n)/\sqrt{n}$
and $g$ be the characteristic function of $G$, the $N(0, \Sigma)$
distribution. From Lemmas
\ref{lm:lm2.1_bp}--\ref{lm:iid_ch}, 
there exist random vectors $S$ and $Y$ on the same probability space
such that $S$ is distributed as $(X_1 + \cdots+ X_n)/\sqrt{n}$,
$Y$ is multivariate normal with mean 0 and variance $\Sigma$ and
\[
P(\| S-Y \| \geq\alpha_1) \leq\alpha_1,
\]
where
\begin{eqnarray*}
\alpha_1 & = & \sigma T + 3(2^{d_1}) e^{-3T^2/32} + \frac{0.25
a_3^3}{\sqrt{n}}  \biggl( \frac{2}{\pi}  \biggr)^{d_1/2} \frac
{T^{d_1}}{\sigma^{d_1+3}}
E(\chi^2(d_1))^{3/2} \\
& &{} + 2 \biggl( \frac{2}{\pi}  \biggr)^{d_1/2}
\frac{T^{d_1}}{\sigma^{d_1}} P  \biggl(\chi^2(d_1) \geq\frac{0.16
n\sigma
^2}{C^2}  \biggr)
+ P\bigl( \| N(0, \Sigma) \| \geq T/2\bigr).
\end{eqnarray*}
From the facts that $E(\chi^2(d_1))^{3/2} \leq(E(\chi
^2(d_1))^2)^{3/4}$ and $P( \| N(0, \Sigma) \| \geq T/2) \leq P (\chi^2(d_1)
\geq T^2/(4a_2^2))$, (\ref{eq:chisq}) and the condition
$a_2\geq1$, we have
\begin{eqnarray*}
\alpha_1 & \leq& \sigma T + 4(2^{d_1}) e^{-3T^2/(32a_2^2)}
+ \frac{0.25a_3^3}{\sqrt{n}}  \biggl( \frac{2}{\pi}  \biggr)^{d_1/2}
\frac{T^{d_1}}{\sigma^{d_1+3}} (2d_1 + d_1^2)^{3/4} \\
& &{} + 2\biggl( \frac{2}{\pi}  \biggr)^{d_1/2} \frac{T^{d_1}}{\sigma
^{d_1}} (2^{d_1}) e^{-0.06n\sigma^2/(C^2)}.
\end{eqnarray*}
Set $\sigma= T^{-1} e^{-3T^2/32}$, then $0< \sigma\leq1$, $T/\sigma
< 12 e^{T^2/8}$ and $1/\sigma< 3e^{T^2/8}$, which, together with
the fact that $( 2/\pi)^{d_1/2} (2d_1 + d_1^2)^{3/4} < 5$, gives that
\begin{eqnarray*}
\alpha_1 &\leq& \bigl(1+4(2^{d_1})\bigr) e^{-3T^2/(32a_2^2)} + \frac{33.75
a_3^3}{\sqrt{n}} (12)^{d_1} e^{(d_1+3)T^2/8} \\
& &{} +2 (19.15)^{d_1} e^{d_1T^2/8} e^{-0.06n\sigma^2/(C^2)} \\
& \leq& \frac{33.75a_3^3}{\sqrt{n}} (12)^{d_1} e^{(d_1+3)T^2/8} +
(48)^{d_1} e^{-3T^2/(32a_2^2)} \leq\alpha,
\end{eqnarray*}
if $0.06n\sigma^2/(C^2) \geq d_1T^2/8 + 3T^2/(32a_2^2)$, which
corresponds to $n \geq(25/(16\times\break a_2^2) + 25d_1/12)C^2T^4
\exp(3T^2/16)$ and we have Lemma \ref{lm:iid_normal}.
\end{pf*}

\subsection[Proof of Theorem 3.1]{Proof of Theorem \protect\ref{thm:consistency}}
\label{sec:proof_thm_consistency}
To prove Theorem \ref{thm:consistency}, we apply Lemma \ref{lm:normal}
by taking the $f_{n,j}(X,Y,z)$'s to be the functions
$\phi^*_{\ell}(X) \phi^*_{\ell'}(X)$, $\phi^*_{\ell}(X) \psi^*_m(Y)$
and $\psi^*_m(Y) \psi^*_{m'}(Y)$, where $1 \leq\ell\leq\ell' \leq p_n$
and $1 \leq m \leq m' \leq q_n$. In such case, (\ref{eq:bdd}) holds
under conditions (B1) and (B2). To see this, for each $1 \leq k \leq n_Z$
and $1 \leq j \leq p_n$, let $\phi^*_{n, j,k}$ be the $j$th component
of $\phi^*$ when $z=z_k$. Then
$\phi^*_{n, j,k} (x) = \sum_{i=1}^{p_n} a_{n, i,j,k} \phi_{n,i}(x)$ for
some $a_{n, i,j,k}$'s and
\begin{eqnarray*}
1 & = &  E  \bigl( (\phi_{n,j,k}^*(X) )^2 |Z=z_k  \bigr) \\
& = & E  \Biggl(  \Biggl( \sum_{i=1}^{p_n} a_{n, i,j,k} \phi_{n,i}(X)
 \Biggr)^2\Bigg|Z=z_k  \Biggr) \\
&\geq& \delta_n \sum_{i=1}^{p_n} a_{n, i,j,k}^2,
\end{eqnarray*}
so $| \phi_{n,j,k}^*(x)| \leq\sqrt{ \sum_{i=1}^{p_n} a_{n, i,j,k}^2 }
\sqrt{\sum_{i=1}^{p_n} \phi_{n,i}^2(x)} \leq\sqrt{p_n/\delta_n}$.
Similarly, for each $1 \leq k \leq n_Z$ and $1 \leq j \leq q_n$, let
$\psi^*_{n, j,k}$ be the $j$th component of $\psi^*$ when $z=z_k$, then
$| \psi_{n,j,k}^*(x)| \leq\sqrt{q_n/\delta_n}$. Thus, (\ref
{eq:bdd}) holds
with $C_n = \max\{ 1, (p_n+q_n)/\delta_n \}$ and it follows from
Lemma \ref{lm:normal} that $\sum_{k=1}^{n_Z} \| \hat{V}^*(z_k) -
V^*(z_k) \|^2$
has the same distribution as $\sum_{k=1}^{n_Z} (n h_n^d c_K
f_Z(z_k))^{-1} \| W_{n,1,k} + W_{n, 2,k} \|^2$, where the $W_{n,1,k}$'s
and $W_{n, 2,k}$'s
are random matrices such that each element in $W_{n,1,k}$ is normal
with mean zero and variance bounded by
$C_n^2=( \max\{ 1, (p_n+q_n)/\delta_n \})^2$, and $\sum_{k=1}^{n_Z}
\|
W_{n, 2,k} \|^2= O_P(\exp(-(\ln n)^{1/9}))$. Therefore,
%
%
\begin{equation} \label{eq:resid} \sum_{k=1}^{n_Z} \| \hat{V}^*(z_k) -
V^*(z_k) \|^2 = O_P ((n h_n^d)^{-1} (\ln n)^{1/8} ).
\end{equation}

To control the difference between $g(\hat{V}^*(z_k), \alpha^*)$ and
$g(V^*(z_k), \alpha^*)$ for $1 \leq k \leq n_Z$,
for a $(p_n+q_n) \times(p_n+q_n)$ matrix $U$, let
%
%
\begin{equation} \label{eq:gstar} g_{i,j}^*(U) =
\cases{
g_{i,j}(U), &\quad if $(i,j) = (1,2)$ or $(2,1)$;\cr
g_{i,j}^{-1}(U), &\quad if $(i,j) = (1,1)$ or $(2,2)$.}
\end{equation}
For $1 \leq k \leq n_Z$, let $\Delta_{i,j,k}=g_{i,j}^*(\hat{V}^*(z_k))
-g_{i,j}^*(V^*(z_k))$ for $1 \leq i, j \leq2$. Then from the fact that
$\| AB \| \leq\| A \| \| B \|$
for two matrices $A$ and $B$, we have
%
%
\begin{eqnarray} \label{eq:consistent3}
&& \| g(\hat{V}^*(z_k), \alpha^*) - g(V^*(z_k), \alpha^*) \|
\nonumber
\\
&&\qquad \leq \prod_{i=1}^2 \prod_{j=1}^2  \bigl( \| g^*_{i,j}(V^*(z_k))
\| +
\| \Delta_{i,j,k} \|  \bigr) - \prod_{i=1}^2 \prod_{j=1}^2 \|
g^*_{i,j}(V^*(z_k)) \| \\
&&\qquad\quad{} +\| g_{1,1}(\hat{V}^*(z_k)) -g_{1,1}(V^*(z_k)) \| \| \alpha^*
(\alpha^*)^T \|.\nonumber
\end{eqnarray}
To control the $\Delta_{1,1,k}$ and $\Delta_{2,2,k}$ in (\ref
{eq:consistent3}), the following result is needed.
\begin{fa} \label{fa:matrix}
Suppose that $A$ is a $p \times p$ invertible matrix and $\Delta= A -
I_p$. Then
$\| A^{-1} -I_p + \Delta\| \leq\| A^{-1} -I_p \| \| \Delta\|$ and
\[
\| A^{-1} -I_p \| \leq\frac{ \| \Delta\| }{ 1-\| \Delta\| }
\qquad\mbox{if $\| \Delta\| < 1$.}
\]
\end{fa}
\begin{pf}
Let $B = A^{-1} -I_p$. Then $B= -
\Delta
- B \Delta$, so $\| B+ \Delta\| = \| B \Delta\| \leq\| B \| \|
\Delta
\|$. Also,
%
%
\begin{equation} \label{eq:p2}
\| B \| \leq\| \Delta\| (1+\| B\| ).
\end{equation}
Apply (\ref{eq:p2}) and we have
\[
\| B \| \leq\frac{ \| \Delta\| }{ 1-\| \Delta\| } \qquad\mbox{if
$\|\Delta\| < 1$.}
\]

Since $\| \alpha^* \|=1$ and for $1 \leq k \leq n_Z$,
$g_{1,1}(V^*(z_k)) = I_{p_n}$, $g_{2,2}(V^*(z_k)) = I_{q_n}$ and
$\| g_{1,2}(V^*(z_k)) \|^2 = \| g_{2,1}(V^*(z_k)) \|^2 \leq(p_n+q_n)
$, from (\ref{eq:consistent3}) and Fact \ref{fa:matrix}, we have
\begin{eqnarray*}
&&
\sum_{k=1}^{n_Z} \| g(\hat{V}^*(z_k), \alpha^*) - g(V^*(z_k),
\alpha^*) \|^2 \\
&&\qquad = O_P \bigl((nh_n^d)^{-1} (\ln n)^{1/8} n_Z^2(p_n+q_n)^3\bigr) \\
&&\qquad =
O_P((nh_n^d)^{-1} (\ln n)^{1/4}),
\end{eqnarray*}
which gives (\ref{eq:consistency1}) since $| \hat{\rho}^2(z_k) -
\rho
_{p_n, q_n}^2(z_k) | \leq\| g(\hat{V}^*(z_k), \alpha^*) - g(V^*(z_k),
\alpha^*) \|$
for $1\leq k \leq n_Z$. (\ref{eq:consistency2}) follows from (\ref
{eq:consistency1}) and the fact that $\sum_{k=1}^{n_Z} ( \hat{f}_Z(z_k)
- f_Z(z_k) )^2$
is $O_P(n_Z(nh_n^d)^{-1})$. The proof of Theorem \ref{thm:consistency}
is complete.
\end{pf}

\subsection[Proof of Theorem 3.2]{Proof of Theorem \protect\ref{thm:distribution}}
\label{sec:proof_thm_distribution}
From Lemma \ref{lm:normal}, the joint distribution of $\hat{V}^*(z_k)\dvtx
1\leq k \leq n_Z$ is the same as that of $V^*(z_k) +
(nh_n^d c_K f_Z(z_k))^{-1/2} (W_{n,1,k} + W_{n,2,k})\dvtx1\leq k \leq
n_Z$, where
%
%
\begin{equation} \label{eq:resid_w2}
\sum_{k=1}^{n_Z} \| W_{n,2,k} \|^2 = O_P(\exp(-(\ln n)^{1/9}))
\end{equation}
and $W_{n,1,k}$'s are independent symmetric normal matrices of mean
zero. To describe the covariance structure of each $W_{n,1,k}$, let
$\phi^* =
( \phi^*_{1}, \ldots, \phi^*_{p_n})^T$, $\psi^* = ( \psi^*_{1},
\ldots,
\psi^*_{q_n})^T$ and let $V_0$ be the $(p_n+q_n)\times(p_n+q_n)$
symmetric matrix
such that $g_{1,1}(V_0) = \phi^*(X) \phi^*(X)^T$, $g_{1,2}(V_0) =
\phi
^*(X) \psi^*(Y)^T$ and $g_{2,2}(V_0) = \psi^*(Y)\psi^*(Y)^T$. For
$1\leq k \leq n_Z$
and $1 \leq m, \ell\leq p_n + q_n$, let $U_{k, m, \ell}$ and $V_{0, m,
\ell}$ be the $(m, \ell)$th elements of $W_{n,1,k}$ and $V_0$,
respectively, then
\[
\operatorname{Cov} (U_{k, m, \ell}, U_{k, m', \ell'}) = \operatorname{Cov}(V_{0, m, \ell}, V_{0, m',
\ell'} |Z=z_k)
\]
for $(m, \ell)$, $(m', \ell') \in\{ (i, j)\dvtx 1 \leq i \leq j \leq
(p_n+q_n) \}$.
For $1\leq k \leq n_Z$, let $\tilde{V}_k = V^*(z_k) + (nh_n^d c_K
f_Z(z_k))^{-1/2} (W_{n,1,k} + W_{n,2,k})$ and
\begin{eqnarray*}
A_1(z_k) & = & g(\tilde{V}_k, \alpha^*) g_{1,1}(\tilde{V}_k) \\
& = & g_{1,2}(\tilde{V}_k) (g_{2,2}(\tilde{V}_k))^{-1} g_{2,1}(\tilde
{V}_k)\\
&&{} - g_{1,1}(\tilde{V}_k) \alpha^* (\alpha^*)^T g_{1,1}(\tilde{V}_k),
\end{eqnarray*}
and let $\tilde{\rho}_0^2(z_k)$ be the largest eigenvalue of $A_1(z_k)
(g_{1,1}(\tilde{V}_k))^{-1}$, then the joint distribution of $\hat
{\rho
}^2(z_k)\dvtx
1 \leq k \leq n_Z$ is the same as that of $\tilde{\rho}_0^2(z_k)\dvtx 1
\leq k \leq n_Z$. For $1 \leq i, j \leq2$ and $1 \leq k \leq n_Z$, let
$\Delta_{i,j,k}
= g_{i,j}(\tilde{V}_k) - g_{i,j}(V^*(z_k))$, then from (\ref{eq:resid}),
%
%
\begin{equation} \label{eq:resid_delta}
{\sum_{k=1}^{n_Z} \sum_{i=1}^2 \sum_{j=1}^2} \| \Delta_{i,j,k} \|^2 =
O_P( (n h_n^d)^{-1} (\ln n)^{1/8})
\end{equation}
and
%
%
\begin{eqnarray}\label{eq:a1_1}
A_1(z_k) &=& g_{1,2}(V^*(z_k)) (g_{2,2}(\tilde{V}_k))^{-1}
g_{2,1}(V^*(z_k))\nonumber\\
&&{} - g_{1,1}(\tilde{V}_k) \alpha^* (\alpha^*)^T
g_{1,1}(\tilde{V}_k) + g_{1,2}(V^*(z_k)) \Delta_{2,1,k}\nonumber\\
&&{}
+ \Delta_{1,2,k} g_{2,1}(V^*(z_k))+ \Delta_{1,2,k} \Delta_{2,1,k}\\
&&{}
- g_{1,2}(V^*(z_k)) \Delta_{2,2,k} \Delta_{2,1,k}\nonumber\\
&&{} - \Delta_{1,2,k}
\Delta_{2,2,k} g_{2,1}(V^*(z_k)) + R_{1,n,k},\nonumber
\end{eqnarray}
where
\begin{eqnarray*}
R_{1,n,k} &=& \Delta_{1,2,k} \bigl(g_{2,2}(\tilde{V}_k)^{-1} -I_{q_n}\bigr)
\Delta_{2,1,k} \\
&&{}
+ g_{1,2}(V^*(z_k)) \bigl(g_{2,2}(\tilde{V}_k)^{-1} -I_{q_n} + \Delta
_{2,2,k}\bigr) \Delta_{2,1,k} \\
&&{}
+ \Delta_{1,2,k} \bigl(g_{2,2}(\tilde{V}_k)^{-1} -I_{q_n}+ \Delta_{2,2,k}\bigr)
g_{2,1}(V^*(z_k)).
\end{eqnarray*}
To simplify the expression for $A_1(z_k)$ in (\ref{eq:a1_1}), we will
make use of the following properties.
\begin{enumerate}[(C6)]
\item[(C4)] The elements of the matrix $g_{1,2}(V^*(z_k))$ are zeros
except that the $(1,1)$th element is 1.
\item[(C5)] For $(i,j) \in\{ (1,2), (2,1) \}$, $g_{i,j}(V^*(z_k))$'s
first row (or first column) is either the first row or the first column
of $g_{i',j'}(V^*(z_k))$ for $(i', j') \neq(i,j)$.
\item[(C6)] The $(1,1)$th element in $g_{2,2}(\hat{V}^*(z_k))$ is 1.
\end{enumerate}
Here (C4) follows from the conditional independence assumption and
(\ref{eq:new_basis_3}), and
(C5) and (C6) follow from (\ref{eq:new_basis_1}). From (C6),
$g_{2,2}(\tilde{V}_k)$ can be expressed as
\[
g_{2,2}(\tilde{V}_k) =  \pmatrix{
1 & B_k^T \cr B_k & D_k}
\]
for some matrices $B_k$ and $D_k$, so the $(1,1)$th element of
$g_{2,2}(\tilde{V}_k)^{-1}$ is $(1 + B_k^T(D_k-B_kB_k^T)^{-1} B_k)$.
Let $J = \alpha^* (\alpha^*)^T$, then by (C4) and (C5), we have
\[
g_{1,2}(V^*(z_k)) (g_{2,2}(\tilde{V}_k))^{-1} g_{2,1}(V^*(z_k)) = \bigl(1 +
B_k^T(D_k-B_kB_k^T)^{-1} B_k\bigr) J,
\]
$g_{1,2}(V^*(z_k)) \Delta_{2,1,k} =J \Delta_{1,1,k}$ and $B_k^TB_k J =
g_{1,2}(V^*(z_k)) (\Delta_{2,2,k})^2g_{2,1}(V^*(z_k))$,
so the expression for $A_1(z_k)$ in (\ref{eq:a1_1}) becomes
\begin{eqnarray*}
&& B_k^T\bigl((D_k-B_kB_k^T)^{-1} - I_{q_n-1}\bigr)B_k J+ g_{1,2}(V^*(z_k))
(\Delta_{2,2,k})^2 g_{2,1}(V^*(z_k)) \\
&&\qquad{} - \Delta_{1,1,k} g_{1,2}(V^*(z_k)) g_{2,1}(V^*(z_k)) \Delta_{1,1,k}
+ \Delta_{1,2,k} \Delta_{2,1,k} \\
&&\qquad{}
- g_{1,2}(V^*(z_k)) \Delta_{2,2,k} \Delta_{2,1,k} - \Delta_{1,2,k}
\Delta_{2,2,k} g_{2,1}(V^*(z_k)) + R_{1,n,k}.
\end{eqnarray*}
Let
\begin{eqnarray*}
A_2(z_k) & = & g_{1,2}(V^*(z_k)) (g_{2,2}(W_{1,n,k}))^2
g_{2,1}(V^*(z_k)) \\
&&{} - g_{1,1}(W_{1,n,k}) g_{1,2}(V^*(z_k)) g_{2,1}(V^*(z_k))
g_{1,1}(W_{1,n,k})\\
&&{} +g_{1,2}(W_{1,n,k}) g_{2,1}(W_{1,n,k})
- g_{1,2}(V^*(z_k)) g_{2,2}(W_{1,n,k}) g_{2,1}(W_{1,n,k})\\
&&{} -
g_{1,2}(W_{1,n,k}) g_{2,2}(W_{1,n,k}) g_{2,1}(V^*(z_k))
\end{eqnarray*}
and
\begin{eqnarray*}
R_{2,n,k} & = & B_k^T\bigl((D_k -B_k B_k^T)^{-1} - I_{q_n-1}\bigr) B_k J \\
&&{} - (n h_n^dc_K f_Z(z_k))^{-1} A_2(z_k) + g_{1,2}(V^*(z_k)) (\Delta
_{2,2,k})^2 g_{2,1}(V^*(z_k))
\\
&&{} - \Delta_{1,1,k}g_{1,2}(V^*(z_k)) g_{2,1}(V^*(z_k)) \Delta
_{1,1,k} + \Delta_{1,2,k} \Delta_{2,1,k} \\
&&{} - g_{1,2}(V^*(z_k)) \Delta_{2,2,k} \Delta_{2,1,k} - \Delta_{2,1,k}
\Delta_{2,2,k} g_{2,1}(V^*(z_k)),
\end{eqnarray*}
then
%
%
\begin{equation} \label{eq:a1} A_1(z_k) = \frac{A_2(z_k)}{n h_n^dc_K
f_Z(z_k)} + R_{1,n,k} + R_{2,n,k},
\end{equation}
%
where
%
%
\begin{equation} \label{eq:resid_r12}
\sum_{k=1}^{n_Z} ( \| R_{1,n,k} \|^2 + \| R_{2,n, k} \|^2)= O_P \biggl(
\frac{\exp(-(\ln n)^{1/9}) (\ln n)^{1/8} }{(n h_n^d)^2}  \biggr)
\end{equation}
from Fact \ref{fa:matrix}, (\ref{eq:resid_w2}) and (\ref
{eq:resid_delta}), and
a simple expression for $A_2 (z_k)$ can be obtained as stated below in
(C7), which follows from (C4) and (C5).
\begin{enumerate}[(C7)]
\item[(C7)] For $1 \leq k \leq n_Z$, $A_2 (z_k) = C_kC_k^T$, where
$C_k$ is the $p_n \times q_n$ matrix obtained by replacing elements in
the first row and first column of $g_{1,2}(W_{1,n,k})$ with zeros.
\end{enumerate}

Note that from (C7), we have that
\[
{\sum_{k=1}^{n_Z}} \| A_2(z_k) \|^2 = O_P\bigl( n_Z (p_n-1)^2(q_n-1)^2 \bigr) =
O_P((\ln n)^{1/8}),
\]
which, together with (\ref{eq:a1}) and (\ref{eq:resid_r12}), implies
that
%
%
\begin{equation} \label{eq:a1_2} {\sum_{k=1}^{n_Z}} \| A_1(z_k) \|^2
= O_P((nh_n^d)^{-2} (\ln n)^{1/8} ),
\end{equation}
and then it follows from (\ref{eq:a1_2}), Fact \ref{fa:matrix} and
(\ref
{eq:resid_delta}) that
%
%
\begin{equation} \label{eq:a1_3} {\sum_{k=1}^{n_Z} }\|
A_1(z_k)(g_{1,1}(\tilde{V}_k))^{-1} - A_1(z_k) \|^2 = O_p((n
h_n^d)^{-3} (\ln n)^{1/4}).
\end{equation}
For $1 \leq k \leq n_Z$, let $\lambda_{0,k}$ be the largest eigenvalue
of $A_2 (z_k)$ and recall that
$\tilde{\rho}_0^2(z_k)$ is the largest eigenvalue of
$A_1(z_k)(g_{1,1}(\tilde{V}_k))^{-1}$. Then by (\ref{eq:a1}), (\ref
{eq:resid_r12}) and (\ref{eq:a1_3}),
%
%
\begin{equation} \label{eq:dist}\qquad
\sum_{k=1}^{n_Z} \bigl(n h_n^d c_K f_Z(z_k) \tilde{\rho}_0^2(z_k) -
\lambda
_{0,k}\bigr)^2 = O_P ( \exp(-(\ln n)^{1/9}) (\ln n)^{1/8}  ).
\end{equation}
Let $\tilde{f}_k$, $\tilde{\rho}(z_k)$ and $\lambda_k\dvtx1 \leq k
\leq
n_Z$ be random variables such that
the joint distribution of $(\tilde{f}_k, \tilde{\rho}(z_k))\dvtx1
\leq k
\leq n_Z$ is the same as that of $(\hat{f}_Z(z_k), \hat{\rho}(z_k))\dvtx
1 \leq k \leq n_Z$,
and the joint distribution of $(\tilde{\rho}(z_k), \lambda_k)\dvtx1
\leq
k \leq n_Z$ is the same as that of $(\tilde{\rho}_0(z_k), \lambda_{0,k})\dvtx
1 \leq k \leq n_Z$. Note that from (\ref{eq:dist}) and the fact that
\[
\sum_{k=1}^{n_Z} \| A_2(z_k) \|^2 = O_P\bigl( n_Z (p_n-1)^2(q_n-1)^2 \bigr),
\]
we have that
\[
\sum_{k=1}^{n_Z} n h_n^d c_K f_Z(z_k) \tilde{\rho}^2(z_k) = \sqrt{ O_P\bigl(
n_Z^2 (p_n-1)^2(q_n-1)^2 \bigr) } = O_P((\ln n)^{1/16}),
\]
so
$n h_n^d \sum_{k=1}^{n_Z} (\hat{\rho}(z_k))^2 = O_P((\ln n)^{1/16})$,
\begin{eqnarray*}
&&  \Biggl| nh_n^d c_K \sum_{k=1}^{n_Z} \hat{f}_Z(z_k) (\hat{\rho
}(z_k))^2 - nh_n^d c_K \sum_{k=1}^{n_Z} f_Z(z_k) (\hat{\rho}(z_k))^2
\Biggr| \\
&&\qquad \leq nh_n^d c_K  \Biggl( \sum_{k=1}^{n_Z}  \bigl( \hat{f}_Z(z_k) -
f_Z(z_k)  \bigr)^2  \Biggr)^{1/2} \sum_{k=1}^{n_Z} (\hat{\rho
}(z_k))^2 \\
&&\qquad = O_P((\ln n)^{1/16})  ( O_P(n_Z (nh_n^d)^{-1})  )^{1/2}\\
&&\qquad = O_P((nh_n^d)^{-1/2} (\ln n)^{3/32})
\end{eqnarray*}
and
\begin{eqnarray*}
&& \hspace*{-22pt} \Biggl| nh_n^d c_K \sum_{k=1}^{n_Z} \tilde{f}_k (\tilde{\rho
}(z_k))^2 - \sum_{k=1}^{n_Z} \lambda_k  \Biggr| \\
&& \leq O_P((nh_n^d)^{-1/2} (\ln n)^{3/32}) +  \Biggl| nh_n^d c_K \sum
_{k=1}^{n_Z} f_Z(z_k) (\tilde{\rho}(z_k))^2 - \sum_{k=1}^{n_Z}
\lambda
_k  \Biggr| \\
{[\mbox{by (\ref{eq:dist})}]} && \leq O_P((nh_n^d)^{-1/2} (\ln
n)^{3/32}) + \sqrt{ n_Z}  ( O_P ( \exp(-(\ln n)^{1/9}) (\ln
n)^{1/8}  )  )^{1/2} \\
&& = O_P ( \exp(-0.5 (\ln n)^{1/9}) (\ln n)^{3/32}  ).
\end{eqnarray*}
The proof of Theorem \ref{thm:distribution} is complete.

\subsection[Proof of Corollary 1]{Proof of Corollary \protect\ref{cor:normal}}
\label{sec:proof_cor_normal}
To prove Corollary \ref{cor:normal}, it is sufficient to establish
(\ref
{eq:cor_normal_var}) and (\ref{eq:lambda_normal}). To see this, let
$\tilde{f}_k$, $\tilde{\rho}^2(z_k)$ and $\lambda_k\dvtx1 \leq k
\leq
n_Z$ be as in Theorem~\ref{thm:distribution}, then
\[
\frac{ nh_n^d c_K \sum_{k=1}^{n_Z} \hat{f}_Z(z_k) \hat{\rho}^2
(z_k) -
n_Z \mu_{p_n,q_n} }{ \sqrt{ n_Z \sigma^2_{p_n, q_n}} }
\]
has the same distribution as
\begin{eqnarray*}
& & \frac{ nh_n^d c_K \sum_{k=1}^{n_Z} \tilde{f}_k \tilde{\rho}^2(z_k)
- n_Z \mu_{p_n,q_n} }{ \sqrt{ n_Z \sigma^2_{p_n, q_n}} } \\
&&\qquad = {\underbrace{ \frac{ nh_n^d c_K \sum_{k=1}^{n_Z} \tilde{f}_k
\tilde{\rho}^2(z_k) - \sum_{k=1}^{n_Z} \lambda_k}{ \sqrt{ n_Z
\sigma
^2_{p_n, q_n}} } }_{I}} + {\underbrace{ \frac{\sum_{k=1}^{n_Z} \lambda_k
- n_Z \mu_{p_n,q_n} }{ \sqrt{ n_Z \sigma^2_{p_n, q_n}} } }_{\mathit{II}}}.
\end{eqnarray*}
Suppose that (\ref{eq:cor_normal_var}) holds, then $I \rightarrow0$ almost
surely by (\ref{eq:nz}) and Theorem \ref{thm:distribution}. Also,
(\ref
{eq:lambda_normal}) says that $\mathit{II}$ converges to $N(0,1)$ in
distribution. Therefore, (\ref{eq:stat_normal}) holds if (\ref
{eq:cor_normal_var}) and (\ref{eq:lambda_normal}) hold.

To establish (\ref{eq:lambda_normal}), we will verify the Lyapounov
condition
%
%
\begin{equation} \label{eq:lya} \lim_{n\rightarrow\infty} \sum_{k=1}^{n_Z}
\frac{E |\lambda_k - \mu_{p_n, q_n}|^3}{ (n_Z \sigma^2_{p_n,
q_n})^{3/2} } = 0,
\end{equation}
and then apply Lindeberg's central limit theorem. Let $\lambda$ be the
largest eigenvalue of $CC^T$. Then $\lambda\leq \operatorname{tr}(CC^T)$, where
$\operatorname{tr}(CC^T)$ is the trace of $CC^T$, which follows the $\chi^2$
distribution with degrees of freedom $m_{1,n}= (p_n-1)(q_n-1)$. Therefore,
\[
E\lambda^3 \leq E(\operatorname{tr}(CC^T))^3 = m_{1,n} (m_{1,n} +2) (m_{1,n} +4),
\]
which implies that $E|\lambda_1 - \mu_{p_n, q_n}|^3 = O( p_n^3 q_n^3)$,
so (\ref{eq:lya}) follows from (\ref{eq:cor_normal_var}) and (\ref
{eq:lambda_normal}) holds.

It remains to prove (\ref{eq:cor_normal_var}). Consider first the case
where (i) holds. By Theorem~1.1 in Johnstone \cite{MR1863961}, 
%
%
\begin{equation} \label{eq:tracy}
\frac{\lambda_1 - \mu_n }{ \sigma_n } \mbox{ converges in distribution}
\qquad\mbox{as $n \rightarrow\infty$,}
\end{equation}
where
\[
\mu_n = \bigl( \sqrt{q_n-2} + \sqrt{p_n-1} \bigr)^2
\]
and
\[
\sigma_n = \bigl(\sqrt{q_n-2} + \sqrt{p_n-1}\bigr)  \biggl( \frac{1}{q_n-2} +
\frac
{1}{p_n-1}  \biggr)^{1/3}.
\]
Here the limiting distribution is the Tracy--Widom law of order 1. Let
$F$ denote its cumulative distribution function. Suppose that $\epsilon
$, $t_1$ and $t_2$ are real numbers such that $t_1 < t_1 + \epsilon<
t_2 - \epsilon$, which implies that $F(t_2) > F(t _2- \epsilon)$ and
$F(t_1+\epsilon) > F(t_1)$. From (\ref{eq:tracy}),
\[
P  \bigl( \lambda_1 > \mu_n + (t_2 - \epsilon) \sigma_n  \bigr)
\geq1 - F(t_2)
\]
and
\[
P \bigl( \lambda_1 < \mu_n + (t_1 + \epsilon) \sigma_n  \bigr)
\geq F(t_1),
\]
if $n$ is large enough. For such $n$, we have
\[
\sigma^2_{p_n, q_n} \geq\frac{\min(F(t_1), 1-F(t_2))
(t_2-t_1-2\epsilon
)^2 \sigma_n^2}{4},
\]
which gives (\ref{eq:cor_normal_var}). The proof of (\ref
{eq:cor_normal_var}) for the case where (ii) holds
can be done by reversing the roles of $p_n$ and $q_n$. The proof of
Corollary \ref{cor:normal} is complete.

\section*{Acknowledgments}
The author thank Dr. Su-Yun Huang, Dr.
I-Ping Tu and Dr. Hung Chen for helpful discussions, and thank the
reviewers for
constructive comments.

\printaddresses

\end{document}